\documentstyle[12pt]{article}
\newfont{\g}{eufm10}

\newcommand{\gtg}{\mbox{\g g}}

\newcommand{\gtb}{\mbox{\g b}}
\newcommand{\gtp}{\mbox{\g p}}

\newcommand{\hgtg}{\mbox{$\hat{\gtg}$}}

\newcommand{\gta}{\mbox{\g a}}
\newcommand{\gtsl}{\mbox{\g sl}}
\newcommand{\gtosp}{\mbox{\g osp}}

\newcommand{\gtn}{\mbox{\g n}}

\newcommand{\gth}{\mbox{\g h}}
\newcommand{\nc}{\mbox{${\bf C}$}}
\newcommand{\nq}{\mbox{${\bf Q}$}}

\newcommand{\nz} {\mbox{${\bf Z}$}}

\newcommand{\cp} {\mbox{${\bf CP^{1}}$}}

\newcommand{\co}{\mbox{${\cal O}$}}
\newcommand{\tco}{\mbox{$\tilde{\co}$}}

\newcommand{\cd}{\mbox{${\cal D}$}}

\newcommand{\cc}{\mbox{${\cal C}$}}
\newcommand{\ce}{\mbox{${\cal E}$}}
\newcommand{\cf}{\mbox{${\cal F}$}}

\newcommand{\ca}{\mbox{${\cal A}$}}

\newtheorem{theorem}{Theorem}[subsection]
\newtheorem{proposition}[theorem]{Proposition}

\newtheorem{remark}[theorem]{Remark}

\newtheorem{lemma}[theorem]{Lemma}
\newtheorem{corollary}[theorem]{Corollary}

\newtheorem{conjecture}[theorem]{Conjecture}

\title{{ \bf Modular Functor and
Representation Theory of $\widehat{\gtsl_{2}}$ at a Rational Level   }}
\author{Boris Feigin\\
Landau Institute for Theoretical Physics
 \and
Feodor Malikov\\
 Department of Mathematics, University of Southern California }
\date{ }
\begin{document}
\baselineskip 20pt

\maketitle

\begin{abstract}
We define a new modular functor based on Kac-Wakimoto admissible
representations and
the corresponding $\cd-$module on the moduli space of rank 2 vector bundles
with the parabolic
structure. A new fusion functor arises which is related to representation
theory
of the pair ``$\gtosp(1|2),\gtsl_{2}$'' in the same way as the fusion functor
for the Virasoro algebra is related to representation theory of the pair
``$\gtsl_{2},\gtsl_{2}$''.

\vspace{5 mm}

\end{abstract}

\section{Introduction}
In this paper we define a new modular functor based on Kac-Wakimoto admissible
representations
over $\widehat{\gtsl}_{2}$.
The modular functor introduced by Segal \cite{seg} assigns a finite-dimensional
vector space
to the data consisting of a punctured curve, a rank 2 vector bundle and a
collection of integral
dominant highest weights attached to the punctures. Our modular functor does
the same for
the Segal's data (with  integral
dominant highest weights replaced with admissible highest weights) extended by
the lines in the
fibers over the punctures. As the data `` surface, vector bundle, punctures,
lines in fibers over punctures'' evolve, so does the corresponding finite
dimensional vector space.
This leads to a new $\cd-$module on the moduli space of rank 2 vector bundles
with  parabolic
structure (fixed lines in certain fibers). The main feature of this
$\cd-$module, as opposed
to the standard one (see Tsuchiya-Ueno-Yamada \cite{ts_u_ya},
 or Beilinson-Feigin-Mazur \cite{beil_feig_maz}), or Moore-Seiberg
\cite{moorseib}) is that
 it is singular over a certain set
of exceptional vector bundles. The latter is closely related to the Hitchin's
global nilpotent
cone.

We also prove that our $D-$module has (in a proper sense) regular singularities
at infinity and that dimension of the generic fiber can be calculated by the
usual combinatorial algorithm: by pinching the surface the problem is reduced
to the case of a  sphere with $\geq 3$ punctures and further to a collection of
spheres with 3 punctures. Dimension of the space attached to the datum ``3
modules sitting at 3 points on a sphere'' is calculated explicitly. It is a
pure
linear algebra calculation of dimension of the space of coinvariants of a
certain infnite dimensional algebra with coefficients in a certain infinite
dimensional representation. As the  result is amusing we will record it here.

First of all, and it is important, in the genus zero case, one can work with
modules at a generic level, as opposed to admissible representations which only
exist when the  level is rational. It is in complete analogy with the usual
WZW model, where the famous theory of Knizhnik-Zamolodchikov equations arises
from a collection of the so-called Weyl modules sitting on a sphere
(terminology is borrowed from  \cite{kazh_luszt}). The family of Weyl modules
is good for the purpose of studying integrable representations because each
integrable representation is a quotient of some Weyl module.
This is no longer the case as far as admissible representations are concerned.
A family of modules suitable for our needs is that of what we call {\em
generalized Weyl modules}; the latter is defined to be  a Verma module
 quotiented out by a singular vector.

Generalized Weyl modules are naturally parametrized by the symbols
$(V_{r}^{\epsilon},V_{s})\;\; r,s\geq 0,\;\epsilon\in\nz/2\nz$. Here $V_{r}$ is
to be thought of as the $r+1-$dimensional irreducible $\gtsl_{2}-$module;
meaning
of $V_{r}^{\epsilon}$ will be explained soon. It is appropriate to keep in mind
that the conventional Weyl module is defined to be the module induced from
$V_{r}$. Therefore usually Weyl modules are labelled by $\gtsl_{2}-$modules. In
our
situation Weyl modules are those related to symbols $(V_{0}^{0},V_{s})$.

According to Verlinde, dimensions of the spaces associated to 3 modules on a
sphere are structure constants of Verlinde algebra. Result of calculation of
Verlinde algebra in our situation is as follows:

\begin{eqnarray}
\label{intr_ourforgen}
(V_{r_{1}}^{\alpha},V_{s_{1}})\circ(V_{r_{2}}^{\beta},V_{s_{2}})&=&\\
(V_{r_{1}+r_{2}}^{\alpha+\beta},V_{s_{1}}\otimes V_{s_{2}})&+&
(V_{r_{1}+r_{2}-1}^{\alpha+\beta+1},V_{s_{1}}\otimes V_{s_{2}})+
(V_{r_{1}+r_{2}-2}^{\alpha+\beta},V_{s_{1}}\otimes V_{s_{2}})+\cdots
+\nonumber\\
(V_{|r_{1}-r_{2}|}^{\alpha+\beta},V_{s_{1}}\otimes V_{s_{2}})&.&\nonumber
\end{eqnarray}
Recall that the usual Verlinde algebra built on Weyl modules is as follows:

\[V_{s_{1}}\circ V_{s_{2}}=V_{s_{1}}\otimes V_{s_{2}},\]
i.e. it is the Grothendieck ring if the category of finite dimensional
representations of $\gtsl_{2}$. Observe that our formula agrees with the latter
one on Weyl modules.

The first component of the right hand side of our formula is  equally easy to
interpret. It is known that the symbols $V_{r}^{\epsilon}$ naturally
parametrize
finite dimensional representations of the simplest rank 1 superalgebra
$\gtosp(1|2)$. The category of finite dimensional $\gtosp(1|2)-$modules
is a tensor category and (\ref{intr_ourforgen}) reads as follows: Verlinde
algebra is isomorphic to the product of Grothendieck rings of the categories of
finite dimensional representations of $\gtosp(1|2)$ and $\gtsl_{2}$.

It is known in principle what to do when passing from modules  to their
quotients, in our case
from generalized Weyl modules at a generic level to admissible representations
at a rational level: one has to replace Lie algebras with quantized universal
enveloping algebras at roots of unity and consider Grothendieck rings
of the corresponding semi-simple ``quotient categories''. Examples: Verlinde
algebra built on integrable $\gtsl_{2}-$modules has to do with $\gtsl_{2}$ in
this way, and Verlinde algebra built on minimal representations of Virasoro
algebra in this way has to do with 2 copies of $\gtsl_{2}$. It appears that
Verlinde algebra built on admissible representations is related to the pair
$(\gtosp(1|2),\gtsl_{2})$ in exactly the same way as $Vir-$Verlinde algebra
is related to the pair of $\gtsl_{2}$'s.

Interest in admissible representation  originates in the fact that the
characters of admissible representations representations at a fixed level give
a
representation of the modular group. However realization of this fact
immediately gave rise to two puzzles:

(i) Given a representation of the modular group, Verlinde formula  produces
structure constants of Verlinde algebra; in the case of admissible
representations some of the structure constants are negative. This does not
make much sense as they are supposed to count dimensions.

(ii) Quantum Drinfeld-Sokolov reduction provides a functor from the category
of $\widehat{\gtsl}_{2}-$modules to the category of $Vir-$modules, which sends
admissible representations to minimal representations. It should give
an epimorphism (or some weakened version of it) of a suitably defined Verlinde
algebra for  $\widehat{\gtsl}_{2}$ on the well-known Verlinde algebra for
$Vir$.

We are able to give an answer to (ii), and a partial answer to (i).

As far as (ii) is concerned, let us for simplicity step aside and consider
$Vir-$modules at a generic (not necessarily rational) level. Then there is an
analogue of a generalized Weyl module -- Verma module quotiented out by a
singular vector -- and these are naturally parametrized by the symbols
$(V_{r},V_{s})$.  The desired  epimorphism is given by:
\[(V_{r}^{\epsilon},V_{s})\mapsto (V_{r},V_{s})+(V_{r-1},V_{s}).\]

This map is naturally related to the  Drinfeld-Sokolov reduction in the
following way.
As we have fixed the category of representations, we have triangular
decomposition of $\gtsl_{2}$; in particular we have 2 opposite nilpotent
subalgebras, $\nc e$, $\nc f$. Therefore there are in fact 2 Drinfeld-Sokolov
functors, $\phi_{e}$, $\phi_{f}$. It happens that the map above is induced
by the direct sum $\phi_{e}\oplus\phi_{f}$.

As to (i), the situation is as follows. The structure constants naturally
arrange in a tensor $\{c_{ij}^{r}\}$, the indices running through a set of
representations in question. Let us compare the set  $\{c_{ij}^{r}\}$ of the
structure
coefficients of our algebra and the set $\{b_{ij}^{r}\}$ of structure
coefficients of
the algebra calculated by Verlinde formula:

 If our $c_{ij}^{r}=0$, then $b_{ij}^{r}=0$. If $c_{ij}^{r}\neq 0$, then
$b_{ij}^{r}$ is ``most certainly'' zero, however in some exceptional cases it
is
non-zero. The latter cases in our situation are interpreted in the following
way. Recall that we have not only 3 modules, $i,j,r,$ but also 3 Borel
subalgebras, $\gtb_{i},\gtb_{j},\gtb_{r}$, which vary. Now as $c_{ij}^{r}\neq
0$, the fiber of our $D-$module is $\neq 0$ (in fact it is 1-dimensional), if
the
3 Borel subalgebras are pairwise different. If however 2 of them meet, the
fiber
usually vanishes, but sometimes survives. It survives if and only if
$b_{ij}^{r}\neq 0$. If non-zero, $b_{ij}^{r}$ can be $\pm 1$. There is no doubt
that $b_{ij}^{r}$ is a result of some cohomological calculation related to the
$D-$module. Unfortunately we cannot make it more precise at the moment.

Just as in the usual case Weyl modules on a sphere produce a trivial vector
bundle with the flat (Knizhnik-Zamolodchikov)connection,
 in our case
we get a bundle with a flat connection on a space of the 2 times greater
dimension. The extra
coordinates come from the flag manifold, recall that we are dealing  with
moduli of vector bundles
with parabolic structure. Horizontal sections of this connection satisfy a
system of differential equations; we get twice as many equations as
there are KZ equations: half of them are indeed KZ equations and the other half
comes from singular vectors in  Verma modules over $\widehat{\gtsl}_{2}$. The
latter is but natural -- it is exactly one of the lessons of the pioneering
work \cite{bpz}.
 This allows to put the integral formulas for solutions of
Knizhnik-Zamolodchikov equations, which we wrote in \cite{feig_mal}, in a
proper context: they
give horizontal sections of this new connection. We conjecture that our
methods, in fact, provide
all horizontal sections. The relation of our formulas to those in
\cite{sch_varch} is that the
latter are necessarily polynomials as functions on the flag manifold while ours
are not.

We wish to acknowledge that there has been a number of works approaching
WZW model for admissible representation from different points of view, see for
example \cite{awata,feld,f_g_p_p,peters,ramg}. It would be interesting to
relate
 our integral formulas with those in \cite{peters} and the new Hopf algebra of
\cite{ramg} to the above mentioned ``$\gtosp(1|2)\times\gtsl_{2}$'' at roots of
unity. To the best of our knowledge, Verlinde algebras proposed in these work
do not solve (ii) above -- those algebras are rather trivial when compared to
the $Vir-$analogue.  Our starting point, see  \cite{feig_mal1}, was the work
\cite{awata}, where Verlinde algebra for admissible representations was
first calculated (in the form equivalent but much less illuminating than the
one  described above), using
the language which left completely open the problem of existence of a
$D-$module, such that dimension of the  fiber is calculated through this
algebra.

{\bf Acknowledgments.} Parts of this work were reported
at the AMS meeting in Hartford, March, 1995, and at Service de Physique
Theorique at Saclay, in November, 1994. We are grateful to J.-B. Zuber for
invitation and warm hospitality. Considerable part of this work was done over
the 2 years one of us spent at Yale. Inspiring and friendly atmosphere at the
Department of Mathematics contributed a lot -- and so did the discussions
with Igor Frenkel, Ian Grojnowski, Gregg Zuckerman. We are grateful to Itzhak
Bars for bringing to our attention the paper \cite{peters}, to Sanjaye
Ramgoolam
for sending his work, and to David Kazhdan for an interesting conversation at
Harvard.

\section{{\bf Notations and known results}}
\label{notat_known_res}

\subsection{ }
 Some notations from  commutative algebra are as follows:

$\nc [t]$ is a polynomial ring, $\nc [[t]]$ is its completion by
positive powers of $t$; $\nc[t,t^{-1}]$ is a ring of Laurent
polynomials and $\nc ((t))$ is its completion by positive powers of $t$.

By functions on the formal (punctured) neighborhood of a non-singular point on
a curve we
will mean a ring isomorphic to $\nc [[t]]$ ($\nc ((t))$ resp.); to specify
such an isomorphism means to pick a local coordinate $t$. The analogous meaning
will
be given to the phrase `` sections of a vector bundle
on the formal (punctured) neighborhood of a non-singular point on a curve''.

\subsection{ }\label{defofalgebrandmodusu}
Set $\gtg = \gtsl_{2}$, $\hgtg=\widehat{\gtsl}_{2}=\gtsl_{2}
\otimes\nc [z,z^{-1}]\oplus\nc c$. Choose a basis $e,\,h,\,f$ of
$\gtg$
satisfying the standard relations $[h,e]=2e,\;[h,f]=-2f,\;[e,f]=h$.
We say that

$\gtg_{\geq}=\nc e\oplus \nc h$
and
$\hgtg_{\geq} =\gtg\otimes z\nc[[z]]\oplus \gtb\oplus\nc c$ are standard Borel
subalgebras of $\gtg$ and $\hgtg$ resp;

$\gtg_{>}=\nc e$
and
$\hgtg_{>} =\gtg\otimes z\nc[[z]]\oplus \gtg_{>}$ are standard ``maximal
nilpotent subalgebras''
of $\gtg$ and $\hgtg$ resp.;

$\nc h$ and $\nc h\oplus\nc c$ are standard Cartan subalgebras of $\gtg$ and
$\hgtg$ resp.

The Verma module $M_{\lambda,k}$ is a module induced from the
character
of $\gtg\otimes z\nc[[z]]\oplus \gtb\oplus\nc c$ annihilating
$\gtg\otimes z\nc[z]\oplus \nc e$ and sending $ h$ and $c$ to
$\lambda$ and $k$ resp. $k$ is often referred to as
a level. Generator of $M_{\lambda,k}$ is usually denoted by
$v_{\lambda,k}$. A quotient of a Verma module is called highest weight module.

The algebra $\hgtg$ is $\nz^{2}_{+}-$graded by assigning $f\otimes
z^{n}\mapsto
(1,-n),\;e\otimes z^{n}\mapsto (-1,-n)$ and so is a Verma module ( as well as
its quotients):
$M_{\lambda,k}=\oplus
_{i,j}M_{\lambda,k}^{i,j}$.

There is a canonical antiinvolution $\omega :\hgtg\rightarrow\hgtg$
interchanging $\hgtg_{>}$ and
$\hgtg_{<}$ and constant on the Cartan subalgebra. For any highest weight
module $V$ denote by
$V^{c}$ and call {\em contragredient} the module equal to the
restricted dual $V^{\ast}$ as a vector space with the following action of
$\hgtg$:
\[<gx,y>=<x,\omega(g)y>,\;g\in\hgtg,x\in V^{\ast},y\in V.\]

If a highest weight module $V$ is irreducible then it is isomorphic to $V^{c}$.
A morphism of
highest weight modules $V_{1}\rightarrow V_{2}$ naturally induces the morphism
of the
corresponding contragredient modules: $V_{2}^{c}\rightarrow V_{1}^{c}$.

A morphism of Verma modules $M_{\lambda,k}\rightarrow M_{\mu,k}$ is
determined by the image of $v_{\lambda,k}$. The image can be written as
$Sv_{\mu,k}$ for a uniquely determined element $S$ of the universal
enveloping algebra of  $\gtg\otimes z^{-1}\nc[z^{-1}]\oplus\nc f$.
If non-zero, the vector $Sv_{\mu,k}$, or even $S$ for this matter, is called
{\em singular}.
The singular vector can be equivalently defined as an eigenvector of the Cartan
subalgebra
 of $\hgtg$
 annihilated by $\hgtg_{>}$. In this form definition applies to an arbitrary
$\hgtg-$module.

\subsection{Singular vector formula}
\label{Singular_vector_formula}
It follows from Kac-Kazhdan
determinant formula that a singular vector generically appears in the
homogeneous
components
of degree either $n(-1,m),\;m> 0,\,n>0$ or $n(1,m),\;m\geq 0,n> 0$.
Denote
the corresponding singular vectors by $S_{n,m}^{1}$ and
$S_{n,m}^{0}$ resp.

 Singular vectors $S_{nm}^{i}$ were found in ~\cite{malff}
in an unconventional form containing  non-integral powers of
elements
of $\hgtg$ ( see also ~\cite{ba_soch} for another approach):

\begin{equation}
S_{nm}^{1}=(e\otimes z^{-1})^{n+mt}f^{n+(m-1)t}
(e\otimes z^{-1})^{n+(m-2)t}\cdots (e\otimes z^{-1})^{n-mt}
\label{s_v_1},
\end{equation}
\begin{equation}
S_{nm}^{0}=f^{n+mt}(e\otimes z^{-1})^{n+(m-1)t}
f^{n+(m-2)t}\cdots f^{n-mt},
\label{s_v_2}
\end{equation}
where $t=k+2$.

This form is not always convenient to calculate a singular vector. It is,
however,
 a useful tool to derive properties of a singular vector. For example,
denoting
by
$\pi:\hgtg\rightarrow\gtg,\;g\otimes z^{n}\mapsto \gtg$ the evaluation
map,
one uses (~\ref{s_v_1},~\ref{s_v_2}) to derive that (see ~\cite{fuchs}, also
{}~\cite{mal} for
the proof in a more general quantum case):

\begin{eqnarray}
\pi S_{nm}^{1}=(\prod_{i=1}^{m}\prod_{j=1}^{N} P(-it-j))e^{N}
\label{p_s_v_1}\\
\pi S_{nm}^{0}=(\prod_{i=1}^{m}\prod_{j=0}^{N-1} P(it+j))f^{N}
\label{p_s_v_2},
\end{eqnarray}
where $P(t)=ef-(t+1)h-t(t+1)$.

\subsection{Generalized Weyl modules and admissible representations}

The  structure of Verma modules over $\hgtg$ is known in full detail
(\cite{mal_2}).
Outside the critical level ($k=-2$) a Verma module
is generically irreducible. $M_{\lambda,k}$ happens
to be reducible if and only if it contains a singular vector.
If $M_{\lambda,k}$ is reducible then the following 2 cases arise:

(i) $k$ is generic (not rational) and $M_{\lambda,k}$  contains only one
singular vector;

(ii)$k+2=p/q>0$ is a ratio of 2 positive integers and $M_{\lambda,k}$ contains
infinitely
many singular vectors.

It can of course happen that $k+2=p/q<0$. We will not be interested
in this case and confine to mentioning that here the situation is  in a sense
dual to (ii).

\subsubsection{ Case (i)}
\label{Case_(i)}
$M_{\lambda,k}$ contains a unique proper submodule $M$
 generated by the singular vector. $M$ is, in fact, a Verma module.

{\bf Definition.}The irreducible quotient $V_{\lambda,k}$
is called {\em generalized
Weyl module}. $\Box$

\bigskip

There arises the exact sequence
\begin{equation}
\label{exactseqforweyl}
0\rightarrow M\rightarrow M_{\lambda,k}\rightarrow V_{\lambda,k}\rightarrow 0.
\end{equation}

A simple property of Kac-Kazhdan equations \cite{kac_kazhd} is that, given
(\ref{exactseqforweyl}),
 the module
$M$ is irreducible and does not project on any generalized Weyl module. Note
that
if the composition series of a $\hgtg-$module only consist of generalized Weyl
modules
then this module
breaks into a direct sum of its components. (This can be proved by methods of
Deodhar-Gabber-Kac
\cite{deodgabbkac}.)

It is an exercise on Kac-Kazhdan equations to derive that the highest weight
$(\lambda,k)$
 of a generalized Weyl module $V_{\lambda.k}$ belongs to either the line
\begin{equation}
\label{param_eq_line_1}
\lambda=-it+j-1,\;k=t-2,
\end{equation}
for some $i\geq 0,j\geq 1$, or
  to the line
\begin{equation}
\label{param_eq_line_2}
\lambda=it-j-1,\;k=t-2,
\end{equation}
for some $i,j\geq 1$; in both cases $t$ is regarded as a parameter. Formula
(\ref{param_eq_line_1}) cooresponds to the case when $V_{\lambda,k}$ is
obtained from $M_{\lambda,k}$ by quotienting out the singular vector
$S^{0}_{i,j}$;
analogously, (\ref{param_eq_line_2}) cooresponds to the case when
$V_{\lambda,k}$ is obtained from $M_{\lambda,k}$ by quotienting out the
singular vector $S^{1}_{i,j}$.

We see that for a fixed level $k$  generalized Weyl modules are parametrized by
the
triples consisting of a pair of nonnegative numbers, $i$, $j$ in the formulas
above, and an element taking one of the 2 values needed to distinguish between
(\ref{param_eq_line_1}) and (\ref{param_eq_line_2}).
To be more
precise, denote by $V_{i}$ the $i+1-$dimensional irreducible representation of
$\gtg$.

{\bf Notation.}
Assign
to $V_{\lambda,k}$ either the symbol $(V_{i}^{0},V_{j-1}),\; i\geq 0,j\geq 1$
if $(\lambda,k)$ satisfies (\ref{param_eq_line_1}), or the symbol
$(V_{i-1}^{1},V_{j-1}),\; i,j\geq 1$ if $(\lambda,k)$ satisfies
(\ref{param_eq_line_2}). $\Box$

\bigskip

This gives us a one-to-one correspondence between the set of generalized Weyl
modules at a fixed generic level
and the set of symbols $(V_{i}^{\epsilon},V_{j})$, where
$\epsilon$ is understood as an element of $\nz/2\nz$.

 Observe that the conventional Weyl module of the level $k$
 is  defined to be the induced
representation
\[\mbox{Ind}_{\gtg[[z]]\oplus\nc c}^{\hgtg}V_{n},\]
where $\gtg[[z]]$ operates on $V_{n}$ via the evaluation map
$\gtg[[z]]\rightarrow\gtg$
and $c\mapsto k$. From our point of view the Weyl module is a quotient of the
Verma module
$M_{n,k}$ by the submodule generated by the singular vector
$f^{n+1}v_{\lambda,k}$. In other
words, Weyl modules are associated to the symbols $(V_{0}^{0},V_{n})$.
This   partially explains appearance of $\gtg-$modules in our notations.

\subsubsection{Case (ii)}
\label{Case_(ii)}
 A Verma module
contains
infinitely many singular vectors and is embedded in finitely many
other
Verma modules. Among all singular vectors in
 $M_{\lambda,k}$ there are 2 independent ones and these
generate
the maximal proper submodule.
Although formally all such Verma modules look alike a special role is
played
by those which can only embed (non-trivially) in themselves. Highest
weights of such modules were called by Kac and Wakimoto
{\em admissible} (~\cite{kac_wak}) and are described as follows.

Let $k+2=p/q$, where $p,q$ are relatively prime positive integers.
The set of admissible highest weights at the level $k=p/q-2$ is given
by
\[\Lambda_{k}=\{\lambda(m,n)=m\frac{p}{q}-n-1\, :\;0<m\leq q,\,
0\leq n\leq p-1\}.\]

What is said above about the structure of Verma modules implies that
any
Verma module appears in the exact sequence of the form
\begin{equation}
\label{bgg_res}
0\leftarrow L_{\lambda_{0},k}\leftarrow
M_{\lambda_{0},k}\stackrel{d_{0}}{\leftarrow} M_{\lambda_{1},k}\oplus
M_{\mu_{1},k}
\stackrel{d_{1}}{\leftarrow} M_{\lambda_{2},k}\oplus
M_{\mu_{2},k}\stackrel{d_{2}}{\leftarrow}\cdots,
\end{equation}
where $\lambda_{0}$ is an admissible weight at the level $k$ and
$L_{\lambda_{0},k}$ is
the corresponding irreducible module.
$L_{\lambda_{0},k}$ is also called {\em admissible}.
 The exact sequence
(\ref{bgg_res}) is called Bernstein - Gel`fand - Gel'fand ( BGG ) resolution.

Again cohomological arguments show (see e.g. \cite{kac_wak}) that if the
composition series of
a $\hgtg-$module only consists of admissible representations then the module is
completely
reducible.

The parametrization of the set of admissible representations we are going to
use is as follows.
Two different generalized Weyl modules project onto one and the same
admissible representation: formula (\ref{bgg_res}) implies that the the two
modules
projecting onto $L_{\lambda_{0},k}$ are $M_{\lambda_{0},k}/M_{\lambda_{1},k}$
and
$M_{\lambda_{0},k}/M_{\mu_{1},k}$.
Therefore
two different triples $(V_{m}^{\epsilon},V_{n})$
are related to the same admissible represenation.
 Introduce the equivalence relation $\approx$ by
$(V_{m}^{\epsilon},V_{n})\approx (V_{q-1-m}^{\epsilon+1},V_{p-2-n}),\; 0\leq
m\leq q-1,\, 0\leq n\leq p-2$.
 Denote
by $(V_{m}^{\epsilon},V_{n})^{\sim}$ the equivalence class of
$(V_{m}^{\epsilon},V_{n})$.

It easy to check that
admissible representations are parametrized by the equivalence  classes of the
triples:
\begin{equation}
\label{parametrofadmreprformula}
\{\mbox{ admissible representations }\}\Longleftrightarrow
\{(V_{m}^{\epsilon},V_{n})^{\sim}\}.
\end{equation}

\subsection{ }
\label{whatcanbesaidhighrank}
Considerable part of the above carries over to the arbitrary Kac-Moody algebra
case.
Here,
for example, is the definition of an admissible representation.
Drop the condition that $\gtg=\gtsl_{2}$, let $M_{\lambda,k}$
be a Verma module over $\hgtg$ and $L_{\lambda,k}$ be its irreducible quotient.
Call $(\lambda,k)$
admissible if $M_{\lambda,k}$ satisfies the following projectivity condition:
if composition
series of a $\hgtg-$module $W$ contains $L_{\lambda,k}$ then $M_{\lambda,k}$
non-trivially maps
in $W$.

 Unfortunately we do not have a reasonable definition of a generalized Weyl
module
in the higher rank case. This is one of the reasons
for which we have to  confine mostly to the $\gtsl_{2}-$case.

 \subsection{Loop modules }
\label{Loop_modules}

 We will also be using $\hgtg-$modules different from Verma modules
or corresponding irreducible ones.

Denote by $\cf_{\alpha \beta} $ a $\gtg-$ module with the basis
$F_{i},\;i\in\nz$
and the action given by
\[eF_{i}=-(\alpha+i-\beta)F_{i+1},\;hF_{i}=(2\alpha+2i-\beta)F_{i},\;
fF_{i}=(-\alpha-i)F_{i-1}.\]
The space $\cf_{\alpha \beta} ^{\nc^{\ast}}=\cf_{\alpha \beta}\otimes\nc
[z,z^{-1}]$ is endowed with
the natural $\hgtg-$module structure. The elements $F_{ij}=
F_{i}\otimes z^{j},\;i,j\in\nz$
serve as a natural basis in it.

Recall (see \ref{Singular_vector_formula}) that   $S^{1}_{nm},\;S^{0}_{nm}$
stand for a singular
vector of degree $n(-1,m)$ or $n(1,m)$ resp. in a Verma module. The following
formulas are proved
by using (\ref{p_s_v_1},\ref{p_s_v_2}):

\begin{eqnarray}
\label{actonloops1}
S^{1}_{nm}F_{n,nm}=\{&\prod_{i=1}^{m}\prod_{j=1}^{n}&(-it-j-\alpha+\beta)(-it-j-\alpha)\}
\{\prod_{s=1}^{n}
(\alpha+i)\}F_{00}\\
S^{2}_{nm}F_{-n,nm}=\{&\prod_{i=1}^{m}\prod_{j=1}^{n}&(it+j-\alpha+\beta)(-it-j-\alpha)\}
\{\prod_{s=1}^{n}
(\alpha-\beta-i)\}F_{00},\label{actonloops2}
\end{eqnarray}
where $t=k+2$.

\section{Construction of the modular functor.}

Although most of our results have to do with $\gtsl_{2}$,
 up to some point it is no extra effort to work in greater
generality. So until sect.\ref{The_spaces_of_coinvariants}, $\gtg$ will stand
for $\gtsl_{n}$ unless
otherwise stated.
\subsection{ Algebra $\hgtg^{A}$ and categories of $\hgtg^{A}$-modules }
\subsubsection{ }
\label{algebra_hgta}
Let $\cc$ be a smooth compact algebraic curve and $\rho:\ce\rightarrow \cc$ be
a rank
$n$ vector bundle  with a flat connection. The connection relates to a section
$s$ of any
bundle $\ca$ associated with $\ce$ the section $dc$ of $\Omega\otimes\ca$ where
$\Omega$
is the sheaf of differential forms over $\cc$. A typical example of $\ca$ is
the bundle
$\mbox{End}\ce$ of fiberwise endomorphisms of $\ce$. The sheaf of sections of
$\mbox{End}\ce$
is naturally a sheaf of Lie algebras over $\cc$.

For a point $P\in \cc$ let $\gtg^{P}$ be the algebra of sections of
$\mbox{End}\ce$ over the formal
neighborhood of $P$. For a finite subset $\bar{A}=\{ P_{1},P_{2},\ldots ,
P_{m}\}
\subset \cc$ set $\gtg^{\bar{A}}= \oplus_{i=1}^{m}\gtg^{P_{i}}$. Define
$\hgtg^{\bar{A}}$ to
be the central extension of $\gtg^{\bar{A}}$ by the cocycle
 \[ <x,y>=\sum_{i=1}^{m}\mbox{Res}_{P_{i}}\mbox{Tr}dx\cdot y.\]
In particular, we obtain the splitting
\begin{equation}
\hgtg^{\bar{A}}=\gtg^{\bar{A}}\oplus\nc\cdot c.
\label{splitting}
\end{equation}

Consider a finite set $A=\{(P_{1},\gtb_{1}),\ldots ,(P_{m},\gtb_{m})\}$
  where $P_{i}\in\cc$ are pairwise different and $\gtb_{i}$ is a Borel
subalgebra
of the algebra of  traceless linear transformations of the fiber
$\rho^{-1}P_{i}$
 ($1\leq i\leq m$). Let $\bar{A}$ be the projection of $A$ on $\cc$. Set
 $\hgtg^{A}=\hgtg^{\bar{A}}$.

\subsubsection{ }
\label{defofocat}
Given $A$ as above, set $\gtn_{i}=[\gtb_{i},\gtb_{i}]$.
Denote by $\hgtg^{A}_{>}$ the subalgebra
consisting of sections $x(.)$ such that $x(P_{i})\in \gtn_{i},\;1\leq i\leq m,$
and by $\hgtg^{A}_{\geq}$  the subalgebra spanned by the space
 of sections $x(.)$ such that $x(P_{i})\in \gtb_{i},\;1\leq i\leq m,$
 and the central
element $c$.
These are analogues of the maximal ``nilpotent'' and maximal ``solvable''
subalgebras for $\hgtg^{A}$, c.f.\ref{defofalgebrandmodusu}.

Denote by $\co^{A}_{k},\;k\in\nc$, the category of
finitely generated $\hgtg^{A}$-modules satisfying
 the conditions:

(i) $c$ acts as  multiplication by $k$;

(ii) the action of the subalgebra $\hgtg^{A}_{>}$ is locally finite.

In much the same way as in \ref{defofalgebrandmodusu} one defines
Verma and generalized Weyl modules over $\hgtg^{A}$:

{\bf Definition.}

(i) We will say that $(\lambda,k)$ is a highest weight of $\hgtg^{A}$ if
$\lambda$ is a functional on $\oplus _{i}\gtb_{i}/\gtn_{i}$ and $k$ is a
number.

(ii) A highest weight $(\lambda,k)$ naturally determines a character of
$\hgtg^{A}_{\geq}$
sending $c$ to $k$ and
 annihilating $\hgtg^{A}_{>}$. Denote by $\nc_{\lambda,k}$
the corresponding 1-dimensional representation.

(iii)  Define the Verma module
$M^{A}_{\lambda ,k}$ to be the induced representation
\[\mbox{Ind}_{\hgtg^{A}_{\geq}}^{\hgtg^{A}}\nc_{\lambda}. \; \Box\]

\bigskip

There is an isomorphism
\[M^{A}_{\lambda ,k}\approx
 \otimes _{i=1}^{m} M^{P_{i},\gtb_{i}}_{\lambda_{i},k}.\]

Suppose now that each $M^{P_{i},\gtb_{i}}_{\lambda_{i},k}$ has at least one
singular vector. If
$k\in \nc \setminus \nq$ then this
singular vector is unique for each $i$. Quotienting out all of
them one obtains the
 {\em generalized Weyl module} $V_{\lambda,k}^{A}$. As above there is an
isomorphism
\[V_{\lambda,k}^{A}\approx
\otimes _{i=1}^{m} V^{P_{i},\gtb_{i}}_{\lambda_{i},k}.\]

If $k$ is not a rational number then any generalized Weyl module is
irreducible. Denote by $\tco_{k}$ the full subcategory of $\co_{k}$ consisting
of all $\hgtg^{A}-$modules whose composition series consist of generalized
Weyl modules. Again if $k$ is not a rational number then $\tco_{k}$ is
semisimple.

If $k$ is rational then there arises the admissible representation
$L_{\lambda,k}^{A}$ if $(\lambda,k)$ is admissible. If the composition series
of a module $V^{A}$ consists only of admissible representations, then $V^{A}$
is
completely reducible.

\begin{remark}
\label{functorwrtA}
There is a canonical isomorphism $\gtb_{1}/\gtn_{1}\approx\gtb_{2}/\gtn_{2}$
for any 2 Borel subalgebras $\gtb_{1},\gtb_{2}$. Therefore if Borel subalgebras
appearing in $A$ evolve, so does the projectivization of the module
$M_{\lambda,k}^A$
or its quotients.

  \end{remark}

\subsubsection{ }
\label{algebra_gofA}
Let $A$ be as in \ref{algebra_hgta}. Let $\gtg(\cc,A)$ be the Lie algebra
of meromorphic sections of
 $\mbox{End}\ce$  holomorphic outside $\bar{A}$.
 The  maps of restriction to formal neighborhoods give rise to the Lie algebra
morphism
\begin{equation}
\label{globalfunct_local}
\gtg(\cc,A)\rightarrow \gtg^{A}
\end{equation}
The splitting (\ref{splitting}) provides us with the section
$s_{A}:\; \gtg^{A}\rightarrow \hgtg^{A}$. Composition of
(\ref{globalfunct_local}) with $s_{A}$
gives the linear morphism
\begin{equation}
\label{globalfunct_affine}
\gtg(\cc,A)\rightarrow\hgtg^{A}.
\end{equation}
The residue theorem implies that (\ref{globalfunct_affine}) is a Lie algebra
morphism
(even though $s_{A}$ is not!).

By (\ref{globalfunct_affine}), the standard
pullback makes each object of $M^{A}\in\co_{k}^{A}$ into a
$\gtg(\cc,A)$-module.
Hence there arises the space of coinvariants
\[(M^{A})_{\gtg(\cc,A)}=M/\gtg(\cc,A)M.\]

\subsection{ Localization of $\hgtg^{A}$-modules }

\subsubsection{ }
\label{whi1jetparam}
  let us
recall that with an $n-$dimensional vector space $W$ one associates the flag
manifold $F(W)=
GL(n,\nc)/B$ and {\em the base affine space} $Base(W)= GL(n,\nc)/N$, where $B$
is a Borel
subgroup and $N$ unipotent subgroup of $B$. The natural map $Base(W)\rightarrow
F(W)$ is a
principal $(\nc^{\ast})^{\times n}$-bundle.

Now return to a $\hgtg^{A}-$module $V^{A}$ and suppose for simplicity that $A$
consists of
1 element $(P,\gtb)$. Consider a family of the data
$\{P, \ce\rightarrow\cc\}$ -- let us not care about Borel subalgebras for the
moment. One expects that the corresponding family of vector spaces arranges
then in a locally trivial vector bundle. An obstacle to get this is that we
have defined
$V^{P}$ up to an isomorphism but have not specified any such isomorphism. For
example, an attempt to
choose a basis in $V^{P}$ requires to choose
(in particular) a local coordinate $z$ at $P$, such that $z(P)=0$.
Different choices of $z$ are essentially different as the group
$\mbox{Diff}(P)$
of diffeomorphisms of the formal neighborhood of $P$ does not
in general act on $V^{P}$. However the subgroup $\mbox{Diff}(P)_{1}\subset
\mbox{Diff}(P)$
of diffeomorphisms preserving the 1-jet of parameter does act on $V^{P}$. We
see that $V^{P}$,
in fact, depends on the 1-jet of parameter at $P$.

To take care of Borel subalgebras,  let us
recall that with an $n-$dimensional vector space $W$ one associates the flag
manifold $F(W)=
GL(n,\nc)/B$ and {\em the base affine space} $Base(W)= GL(n,\nc)/N$, where $B$
is a Borel
subgroup and $N$ unipotent subgroup of $B$. The natural map $Base(W)\rightarrow
F(W)$ is a
principal $(\nc^{\ast})^{\times n}$-bundle.

Similar arguments applied to $\gtb$ show that

the module $V^{A}=V^{P,\gtb}$ depends on the quadruple $(P,\gtb,j,x)$ such that
$j$ is a 1-jet
of parameter at $P$ and $x\in Base(\nc^{n})$ belongs to the preimage of $\gtb$.

One concludes that we do get a locally trivial vector bundle after pull-back to
the space of pairs ``1-jet of parameter at $P$, element of the maximal torus of
the
Borel group related to $\gtb$''. Let us be more precise now.

\subsubsection{ }
\label{maingenerresults}

Let $\bar{\pi}:\;\cc_{S}\rightarrow S$ be a family of smooth projective curves
and $\rho_{S}:\; \ce_{S}\rightarrow \cc_{S}$ be a rank $n$ vector bundle. There
arise 2 more
bundles:

(i) the  bundle $Base(\rho_{S}):\; Base(\ce_{S})\rightarrow \cc_{S}$ with the
fiber
over any $x\in\cc_{S}$  equal to the base affine space
of the vector space $\rho_{S}^{-1}x$;

(ii) the $\nc^{\ast}-$bundle $J^{(1)}(\cc_{S})\rightarrow \cc_{S}$ of  1-jets
of coordinates along
fibers of $\bar{\pi}$.

 Consider the fibered product
$Base(\ce_{S})\times_{\cc_{S}}J^{(1)}(\cc_{S})$ and the natural map
\[\pi:\; Base(\ce_{S})\times_{\cc_{S}}J^{(1)}(\cc_{S})\rightarrow S.\]

Pick a  non empty  finite set $A_{S}$ of sections of $\pi$   satisfying the
condition:

for any $s\in S$ the natural projection of
 the set $A_{S}(s)=\{a(s),\;a\in A_{S}\}$
 on  $\bar{\pi}^{-1}(s)$
is an injection.

Pick  an arbitrary curve, say $\cc_{s_{0}}$, from our family. Consider a
highest weight
module  $M^{A }$ over $\hgtg^{A}$, where we write $A$ instead of the lengthy
$A_{S}(s_{0})$; what follows is obviously independent of the choice of $s_{0}$.

By \ref{defofocat}, remark  \ref{functorwrtA}, and \ref{whi1jetparam}, we get a
$\hgtg^{A_{S}(s)}-$module $M^{A_{S}(s)}$ for any $s\in S$ and the collection
 $\{M^{A_{S}(s)},\; s\in S\}$ arranges in a locally trivial vector bundle.
With
each $s\in S$ we can further associate a vector space, that is the space of
coinvariants
\[(M^{A_{S}(s)})_{\gtg(\pi^{-1}S,A_{S}(s))},\]
see \ref{algebra_gofA}.

\begin{theorem}
\label{existofDmod}
Suppose the collection $\psi=(M^{A},\pi,A_{S})$, satisfying the conditions
imposed above, is given.
 Then
there is a twisted $ \cd-$module (that is
a sheaf of modules over a certain algebra of twisted differential operators) on
$S$ such
that its fiber over $s\in S$ is $(M^{A_{S}(s)})_{\gtg(\pi^{-1}S,A_{S}(s)}$.
\end{theorem}

This theorem is an immediate consequence of \cite{beil_feig_maz} and
\cite{bern_beil,bryl_kash} . Briefly the construction is as follows. Take a
vector field  $\xi$ on $U\subset S$. It lifts to a meromorphic vector field
on $\cc_{S} - A_{S}(S)$ over $U$, and
further to a meromorphic vector field on $\pi^{-1}(U)\subset
Base(\ce_{S})\times_{\cc_{S}}J^{(1)}(\cc_{S})$; denote this vector field by
$\xi^{\ast}$.  Trivializing  the infinitesimal neighborhood of
$A_{S}(U)\subset\cc_{S}$ by chosing, locally with respect to $U\subset S$,
coordinates in the fibers, one gets vertical components
$\{\xi^{\ast}_{vert;i}\}$, so that $\xi^{\ast}_{vert;i}$ is the vertical
component in the formal neighborhood of the $i-$ section. Projecting
$\xi^{\ast}_{vert;i}$ on $Base(\ce_{S})$ one gets some element of $U(\gtg)$,
say
$u_{i}$;  projecting $\xi^{\ast}_{vert;i}$ on $J^{(1)}(\cc_{S})$ one gets some
vector field, say $v_{i}$. Both $u_{i},\; v_{i}$ act on our $\hgtg^{A}$-module
$M^{A}$:
$u_{i}$ naturally, $v_{i}$ by means of the Sugawara construction. Going over
definitions one gets that this well defines a twisted $\cd-$module with the
fiber as in the theorem. $\Box$

\bigskip

Denote the constructed $\cd-$module by $\Delta_{\psi}(M^{A})$.

In the case when $M^{A}$ is an admissible representation the following result
is valid.

\begin{theorem}
\label{smoothinadm}
If $n=2$ and $M^{A}$ is an admissible $\hgtg^{A}-$module then
$\Delta_{\psi}(M^{A})$ is {\em holonomic}
for almost any vector bundle $\ce_{S}$
(i.e.as a sheaf $\Delta_{\psi}(M^{A})$ is isomorphic to a sheaf of sections
 of a certain finite rank  vector bundle over some open set in $S$ ).
\end{theorem}

{\bf Proof.}

To prove this theorem essentially  means to show that the spaces
 $(M^{A_{S}(s)})_{\gtg(\pi^{-1}S,A_{S}(s)},\; s\in S,$ are all finite
dimensional. That will be done
in \ref{finofcoinvhighgensubs},  Proposition \ref{proofoflisse} in the higher
genus case
and in \ref{finofcoinvsferesubs}, Proposition \ref{fincinvforcpadmweyl} for
$\cp$. . We will also give there a precise meaning
to the phrase ``almost any vector bundle'' in Theorem\ref{smoothinadm}. $\Box$

\bigskip

Results of \ref{Quadratic_degeneration} will show that the standard
combinatorial algorithm can be used to calculate the dimension of the fiber of
our
$\cd-$module using the dimensions of the spaces of coinvariants on a sphere
with 3 punctures.
The latter dimensions will be calculated in \ref{The_rational_leve_case}.

\section{The spaces of coinvariants}
\label{The_spaces_of_coinvariants}
 In this section we will be concerned
with the space of coinvariants $(M^{A})_{\gtg(\cc,A)}$
(or  spaces closely related to it ) in the case when $M^{A}$ is either
a generalized Weyl module or an admissible representation. The standard tool to
get finiteness results
about coinvariants is the notion of
 {\em singular support}.

\subsection{Singular support and coinvariants}
\label{singsuppandcoinvvvvv}

Let $\gta$ be a Lie algebra. Universal enveloping algebra $U\gta$ is filtered
in
 the standard way so that the associated graded algebra is $S\gta$. One says
that a filtration of a finitely generated $\gta-$module $V$ is good if (i) it
is
compatible with the filtration of $U\gta$, and (ii) the associated graded
module $Gr\, V$ is finitely generated as an $S\gta$-module.

{\bf Definition} Singular support, $SS V$, of $V$ is the zero set of the
vanishing ideal of the $S\gta$-module $Gr\, V$. $\Box$

\bigskip

Obviously, $SS V$ is a conical subset of $\gta^{\ast}$.

For a subalgebra $\gtn\subset \gta$, call  $V$ an $(\gta,\gtn)$-module if it is
an $\gta-$module and
$\gtn$ acts on $V$ locally nilpotently. Typical example: any module from the
$\co-$category is a $(\hgtg,\hgtg_{>})-$module.

\begin{lemma} (\mbox{ see \cite{beil_feig_maz}})
\label{singsupp-finitenessofcoinv}
Let $\gta$ be a Lie algebra and $\gtp\subset\gta$ be its subalgebra. Denote by
$\gtp^{\perp}$ the
annihilator of $\gtp$ in $\gta^{\ast}$
Let $V$ be an $(\gta,\gtn)$-module.
If $SSM\cap\gtp^{\perp}=\{0\}$ and $dim\, \gta/\gtn\oplus\gtp<\infty$ then
$\mbox{dim}M_{\gtp}<\infty$.
\end{lemma}

Recall that from now on $\gtg=\gtsl_{2}$ unless otherwise stated.

\subsection{ Singular support of $\hgtg^{A}-$modules}
\label{singsuppppofhgta}
Observe that there is an involution $\sigma$ of $\hgtg$ sending $f$ to
$e\otimes z^{-1}$ and
$e\otimes z^{-1}$ to $f$, see
\ref{defofalgebrandmodusu}
for notations. There arises the involution, also denoted by $\sigma$, acting
on the algebras $\hgtg^{A}$ and their duals. This involution is not canonical
but we do not have
to care as our considerations here are purely local.

Denote by $\Omega^{A}$ the space of $\gtg-$ valued differential
forms on the formal neighborhoods of the points from $A$. There is a natural
embedding
$\Omega^{A}\hookrightarrow (\gtg^{A})^{\ast}$ (``take the traces and
 then sum up all the residues!'')

We will make use of 2 subspaces of $\Omega^{A}$: $\Omega^{A}_{reg}$ is all
regular forms and
$\Omega^{A}_{nilp}$ is all forms with values in the nilpotent
cone.
\begin{theorem}
\label{theoronsingsupp}

(i) If $M^{A}$ is a generalized Weyl module then
$SSM^{A}=\Omega^{A}_{reg}\cup\sigma\Omega^{A}_{reg}$.

(ii) (E.Frenkel, B.F.) If $M^{A}$ is an admissible representation then
$SSM^{A}=\Omega^{A}_{nilp}\cup\sigma\Omega^{A}_{nilp}$.
\end{theorem}

\begin{remark}
It is easy to see that although $\sigma$ is not determined uniquely the
spaces $\sigma\Omega^{A}_{reg}$, $\sigma\Omega^{A}_{nilp}$ are canonical.
For example  $\sigma\Omega^{A}_{reg}$ is the space of forms such that:

they have at most order 1 pole at $\bar{A}$;

their residue at each $P_{i}\in \bar{A}$ belongs to $\gtn_{i}$;

at each $P_{i}\in\bar{A}$ their constant term  belongs to $\gtb_{i}$.
\end{remark}

\subsection{Finiteness of coinvariants --  the higher genus case}
\label{Finiteness_of_coinvariants-hghergen}

\subsubsection{ Hitchin's theorem.}
\label{hitchinstheorem}
First recall a well-known result of  Hitchin, \cite{hitch}.
With a vector bundle $\ce\rightarrow \cc$ associate the map
\begin{eqnarray}
\label{def-fhitchmap}
H(\ce):\;H^{0}(\cc,\Omega\otimes\mbox{End}\ce)\rightarrow
 \oplus_{i=2}^{n}H^{0}(\cc,\Omega^{\otimes i}),\\
X\mapsto \mbox{Tr}X^{i}\nonumber
\end{eqnarray}
Call a bundle $\ce$ {\em exceptional} if $\mbox{ker}H(\ce)\neq 0$. Obviously
$\mbox{ker}H(\ce)$ is exactly the space of global differential forms
with values in nilpotent endomorphisms of the vector bundle $\ce$.

\begin{theorem} (\mbox{Hitchin \cite{hitch}})
\label{hitchi}
 Zero set of the map (\ref{def-fhitchmap}) is a maximal Lagrangian
submanifold in the cotangent bundle of the moduli space of vector bundles over
$\cc$. In particular,
exceptional
vector bundles form a positive codimension algebraic subset of the moduli space
of vector bundles.
\end{theorem}

For us, importance of Theorem \ref{hitchi} is in  that generically a vector
bundle does not allow
a non-trivial global differential form with coefficients in nilpotent
endomorphisms of the bundle.

\subsubsection { Subtracting lines from rank 2 vector bundles.}
\label{subtrlinesfromvectbundles}

 An analogue of subtracting a point
from a line bundle (or, better to say, from  its divisor) is an operation of
subtracting a line from a rank 2  vector bundle.

To a  rank 2 vector bundle $\ce\rightarrow\cc$ one can associate a  module over
the sheaf of regular functions -- the
sheaf of sections of $\ce$.; denote this sheaf by $Sect(\ce)$. This establishes
a one-to-one correspondence between rank 2 vector bundles and rank 2 locally
free modules over the sheaf of regular functions.
 Now fix a line, $l$, in a fiber of $\ce$ over some point $P\in\cc$. Denote by
$S (l)$ a sheaf such that:

(i) $S(l)|_{U}=Sect(\ce)|_{U}$ if $P$ does not belong to $U$;

(ii)  $S(l)|_{U}, P\in U,$ is the space of meromorphic sections of $\ce$ over
$U$ regular outside $P$, having at most order 1 pole at $P$ and such that their
residue at $P$ belongs to the fixed line $l$.

It is obvious that $S(l)$ is a rank 2 locally free module. Therefore it
defines a rank 2  vector bundle. Denote this vector bundle by $\ce(l)$.
If a collection of lines -- $l_{1},l_{2},...,l_{m}$ -- is subtracted, then
denote
the corresponding vector bundle by $\ce(l_{1}+\cdots + l_{m})$.

Suppose we have a moduli space of rank 2 vector bundles with parabolic
structure with fixed determinant. Elements of such a space are isomorphism
classes of the data
(vector bundle $\ce$, fixed lines $l_{1},...,l_{m}$ in some fibers.) It is
rather
clear  that the map $(\ce,\; l_{1},...,l_{m})\mapsto (\ce(l_{1}+\cdots
+l_{m}),\; l_{1},...,l_{m})$ is a homeomorphism of 2 moduli spaces with
different determinants.

{\bf Definition.}  Call the data $(\ce,\; l_{1},...,l_{m})$ generic if
 $\ce(l_{i_{1}}+\cdots +l_{i_{s}})$ is not  exceptional for
any subset $\{i_{1},...,i_{s}\}\subset\{1,2,..., m\}$. $\Box$

\bigskip

It follows from Theorem \ref{hitchi} that the set of generic vector bundles is
open and everywhere dense.

\subsubsection { Finiteness of coinvariants.}
\label{finofcoinvhighgensubs}

Suppose we are in the situation of \ref{defofocat}: we have an admissible
$\hgtg^{A}-$module
$M^{A}$ on the curve $\cc$ with a vector bundle $\ce\rightarrow\cc$.
As $A$ is a collection of borel subalgebras $\gtb_{1},...,\gtb_{m}$ operating
in fixed fibers, we have
parabolic structure -- lines $l_{1},...,l_{m}$ in the corresponding fibers
preserved by the $\gtb_{i}$'s.
Call the data $(\ce,A)$  generic if the data $(\ce, l_{1},...,l_{m})$ is
generic
in the sense of \ref{subtrlinesfromvectbundles} above.

Recall that we are interested in the space of coinvariants
$M^{A}_{\gtg(\cc,A)}$, where $\gtg(\cc,A)$ is an algebra of endomorphisms of
the bundle $\ce$ regular outside points from the corresponding $\bar{A}$, see
\ref{algebra_gofA} and \ref{maingenerresults}, Theorem \ref{existofDmod}.

\begin{proposition}
\label{proofoflisse}
 Let $(\ce, A)$ be generic. Then
\[\mbox{dim} M^{A}_{\gtg(\cc,A)}<\infty.\]
\end{proposition}

{\bf Proof.}  One extracts from definitions that the annihilator
$\gtg(\cc,A)^{\perp}$ of the algebra $\gtg(\cc,A)$ is the space
$\Omega_{\cc,A}(\ce)$ of global meromorphic $End(\ce)-$valued differential
forms regular outside
$\bar{A}\subset\cc$.

By Theorem \ref{theoronsingsupp}(ii) we get that
$SSM^{A}\cap\gtg(\cc,A)^{\perp}=\Omega_{nilp}(\ce)\cup\sigma\Omega_{nilp}(\ce)$, where
$\Omega _{nilp}(\ce)$ is the space global nilpotent transformations of $\ce$,
and $\sigma$ is the twist introduced in \ref{singsuppppofhgta}.

Genericity condition means that $\Omega_{nilp}(\ce)=0$, see
\ref{hitchinstheorem} and  \ref{subtrlinesfromvectbundles}.

On the other hand it is easy to see that the operation of subtracting a line
generates the twist $\sigma$ on endomorphisms. (In fact one has to compose
subtracting of a line
with a reflection in the fiber, but this does not change the isomorphism class
of the bundle.) Therefore genericity condition also implies that
$\sigma\Omega_{nilp}(\ce)=0$.

Hence we get that $SSM^{A}\cap\gtg(\cc,A)^{\perp}=0$. And as the space
$\hgtg^{A}_{>}+\gtg(\cc,A)$ is  of finite codimension in $\hgtg^{A}$,
application
of Lemma \ref{singsupp-finitenessofcoinv}, see \ref{singsuppandcoinvvvvv},
completes the proof. $\Box$

\bigskip

\bigskip

In order to study quadratic degenerations we will need the following
stronger finiteness result. Along with the set
$A=\{(P_{1},\gtb_{1}),...,(P_{m},\gtb_{m})\}$, consider the set $ A_{2} =\,
\{(P_{m+1},\gtb_{m+1}), (P_{m+2},\gtb_{m+2})\}$ such that the points
$P_{1},...,P_{m+2}\in\cc$ are different. Denote by $\gtg(\cc,A,A_{2})$ the
subalgebra of $\gtg(\cc,A)$ consisting of functions taking values in
$\gtn_{i}=[\gtb_{i},\gtb_{i}]$ at point $P_{i}$, $i=m+1,m+2$.

\begin{proposition}
\label{highgen_finit_afterpionching}

If $(\ce, A\bigsqcup A_{2})$ is generic and $M^{A}$ is admissible, then
\[dim\, (M^{A})_{\gtg(\cc,A,A_{2})}<\infty.\]

\end{proposition}

{\bf Proof.} We are again going to apply Lemma
\ref{singsupp-finitenessofcoinv}. Observe that
$\gtg(\cc,A,A_{2})^{\perp}$ consists of meromorphic forms on $\cc$ with values
in $End(\ce)$,
regular outside $\{P_{1},...,P_{m+2}\}\subset\cc$, having at most order 1 poles
at $P_{m+1},P_{m+2}$, their residues at the latter points lying in $\gtb_{1}$
( $\gtb_{2}$ resp.).

By Theorem \ref{theoronsingsupp}(ii), $\gtg(\cc,A,A_{2})^{\perp}\cap\, SSM^{A}$
consists of forms with values in nilpotent endomorphisms, satisfying the above
listed global conditions. This implies, in particular, that actually residues
of our forms
belong to $\gtn_{m+1},\gtn_{m+2}$ at $P_{m+!}, P_{m+2}$ resp..

Given an element $\omega\in\gtg(\cc,A,A_{2})^{\perp}\cap\, SSM^{A}$, subtract
some lines from $\ce$ so as to make $\omega$ be everywhere regular. Genericity
condition implies then that $\omega=0$, and application of   Lemma
\ref{singsupp-finitenessofcoinv} completes the proof. $\Box$

\bigskip

\subsection{Finiteness of coinvariants -- the case of $\cp$}
\label{Finiteness_of_coinvariants}

\subsubsection{Generic vector bundles on $\cp$}
\label{genervectbumdonsph}

Let $O(n)$ be the degree $n$ line bundle over $\cp$. It is known, e.g.
\cite{pr_seg}, that any rank 2 vector bundle over $\cp$ is a direct sum
$O(r)\oplus O(s)$ for some $r$, $s$.

As there are no moduli, it is hard to speak about generic vector bundles.
Nevertheless we will call $O(r)\oplus O(s)$ {\em exceptional} if $|r-s|>1$.
Here is
a justification.

\begin{lemma}
\label{whyexceptional}

Let $\ce= O(r)\oplus O(s)$ and $(\ce,l_{1},...,l_{m})$, $m\geq |r-s|$, a vector
bundle with parabolic structure. Then generically with respect to
$l_{1},...,l_{m}$ the bundle $\ce(l_{1}+\cdots +l_{m})$ is not exceptional:

\[\ce(l_{1}+\cdots +l_{m})=\left\{\begin{array}{lll}
O(p+1)\oplus O(p)&\mbox{ if } r+s-m=2p+1\\
 O(p)\oplus O(p)&\mbox{ if } r+s-m=2p.
\end{array}
\right.\]
\end{lemma}

Lemma \ref{whyexceptional} seems to be common knowledge, although we failed
to find a reference with its proof.

Proceed just like we did in \ref{subtrlinesfromvectbundles}:  call  $(\ce,\;
l_{1},...,l_{m})$ generic if
 $\ce(l_{i_{1}}+\cdots +l_{i_{s}})$ is not  exceptional for
any subset $\{i_{1},...,i_{s}\}\subset\{1,2,..., m\}$.

\subsubsection{Finiteness of coinvariants}
\label{finofcoinvsferesubs}

A specific feature of the genus zero case is that we do not necessarily have
to consider admissible representations -- generalized Weyl modules, see
\ref{Case_(i)},
 will also do.

Let us again consider a vector bundle $\ce$ over $\cp$ and a $\hgtg^{A}-$module
$M^{A}$. As in \ref{finofcoinvhighgensubs}, $A$ determines a parabolic
structure
on $\ce$, say $(\ce, l_{1},...,l_{m})$. Call the data $(\ce,A)$ generic
if $(\ce, l_{1},...,l_{m})$ is also.

\begin{proposition}
\label{fincinvforcpadmweyl}
If  $(\ce,A)$ is generic and $M^{A}$ is either admissible or generalized Weyl
module, then

\[\mbox{dim } (M^{A})_{\gtg(\cp,A)}<\infty.\]

\end{proposition}

{\bf Proof}  is a simplified version of the proof of Proposition
\ref{proofoflisse} in \ref{finofcoinvhighgensubs}. The new features are as
follows: to include generalized Weyl modules one uses Theorem
\ref{theoronsingsupp}(i) in addition to Theorem \ref{theoronsingsupp}(ii);
instead of the Hitchin's theorem one uses the ``observation'' that $O(n)$ has
no
non-zero global sections if $n<0$. $\Box$

\begin{corollary}
\label{existofdmodforweyl}
 If $M^{A}$ is a generalized Weyl module then
there is a holonomic
twisted $D$-module living in the space $(J^{(1)}(\cp)\times
J^{(1)}(\cp))^{\times m}$ with the fiber $M^{A}_{\gtg(\cp,A)})$.
\end{corollary}

{\bf Proof.} Repeating word for word proof of  Theorem\ref{smoothinadm} one
derives
 from Proposition \ref{fincinvforcpadmweyl} existence of a twisted $D-$module
on
the space $(Base(\nc^{2})\times J^{(1)}(\cp))^{\times m}$. But for $\gtsl_{2}$,
 the flag manifold is $\cp$ and
the base affine space $(Base(\nc^{2})$ is  also the space of 1-jets of
parameter $ J^{(1)}(\cp))^{\times m}$. $\Box$

\bigskip

As in \ref{finofcoinvhighgensubs}, we want to prove a generalization of
Proposition \ref{fincinvforcpadmweyl} in order to prepare grounds for studying
quadratic degeneration.

Along with $A=\{(P_{1},\gtb_{1}),\ldots ,(P_{m},\gtb_{m})\}$ consider 2 sets
$A_{1}=\{(P_{m+1},\gtb_{m+1})\}$ and $A_{2}=\{(P_{m+1},\gtb_{m+1}),
(P_{m+2},\gtb_{m+2})\}$
such that $P_{1},...,P_{m+2}$ are different points in $\cc$.

With $A_{1}$ and $A_{2}$ associate the following 2 subalgebras of
$\gtg(\cp,A)$:
$\gtg(\cp,A,A_{1})$ consists of all functions taking values in
$\gtn_{m+1}=[\gtb_{m+1},\gtb_{m+1}]$
at the point $P_{m+1}$;
$\gtg(\cp,A,A_{2})$ consists of all functions taking values in
$\gtn_{i}=[\gtb_{i},\gtb_{i}]$
at the point $P_{i}$, $i=m+1,\;m+2$.

\begin{proposition}
\label{finofcoinv_general}
 Let   $(\ce, A_{2})$ be
generic. Then

(i) If $M^{A}$ is a generalized Weyl module over $\hgtg^{A}$, then
$\mbox{dim}(M^{A})_{\gtg(\cp,A,A_{1})}<\infty$;

(ii)  If $M^{A}$ is an admissible representation of $\hgtg^{A}$, then
$\mbox{dim}(M^{A})_{\gtg(\cp,A,A_{2})}<\infty$.
\end{proposition}

 \bigskip

{\bf Proof.} of (ii) repeats almost word for word that of Proposition
\ref{highgen_finit_afterpionching} in \ref{finofcoinvhighgensubs} with
simplifications analogous to those
indicated in the proof of Proposition  \ref{fincinvforcpadmweyl}.

As to (i), its proof is again application of the same technique in a slightly
different
form: one has to take a form $\omega\in\gtg(\cp,A,A_{2})^{\perp}\cap SSM^{A}$
and to subtract lines from $\ce$ so as to make $\omega$ into
a form with either one pole (at $P_{m+1}$) or 2 poles (one of them is
again at $P_{m+1}$) in such a way that the bundle obtained is $O(n)\oplus
O(n)$.
The 2 cases are of course distinguished by the parity of the difference
between the degrees of the determinant of $\ce$ and $\omega$. In both cases
it is easy to prove that $\omega=0$ using the fact that any differential form
with trivial coefficients has at least 2 poles. $\Box$

\subsubsection{Holonomic $D$-module on $(\nc\times\nc)^{\times m}$}
\label{holondmodoncc}

 We will now get rid of twisted differential operators in Corollary
\ref{existofdmodforweyl} under the assumption that the vector bundle
$\ce\rightarrow \cp$ is trivial.
Consider the set $A'=A\bigsqcup (P_{\infty},\gtb_{\infty})$. Attach to the
point $(P_{\infty},\gtb_{\infty})$  the module $(V^{0}_{0},V_{0})$ known
as the vacuum representation, see \ref{Case_(i)} for notations.
$(P_{\infty},\gtb_{\infty})$ can be redefined as the module induced from the
trivial representation (see also \ref{Case_(i)}) and therefore there is an
isomorphism $ M^{A}_{\gtg(\cp,A)}\approx M^{A'}_{\gtg(\cp,A')}$. Now consider
the twisted $D-$module with fiber $ M^{A'}_{\gtg(\cp,A')}$ on the space
$(J^{(1)}(\cp)\times J^{(1)}(\cp))^{\times m+1}$. Restrict it to the space
 $(J^{(1)}(\cp)\times J^{(1)}(\cp))^{\times m}$ by having  the point
$(P_{\infty},\gtb_{\infty})$ fixed. The result of this operation is that the
bundles in question trivialize: $\cp - \gtb_{\infty}=\cp - P_{\infty}=\nc$ and
$J^{(1)}(\nc)=\nc^{\ast}\times \nc$. Further pushing forward by ``integrating
along $\nc^{\ast}$'' one gets a $D-$module over the space $(\nc\times
\nc)^{m}$.
Observing that it is appearance of the bundle $J^{(1)}(\cp)\rightarrow\cp$
which
was responsible for the twisting of the $D-$module, one argues that we get
a usual holonomic $D$-module on $(\nc\times \nc)^{m}$ with fiber $
M^{A}_{\gtg(\cp,A)}$. In particular, we get a bundle with flat connection
over an open subset of $(\nc\times \nc)^{m}$.

{\bf Notation.} Denote the constructed in this way bundle with flat connection
by $\Delta(M^{A})$. $\Box$

\bigskip

  We are unable to describe this open
subset explicitly at present. It follows from the requirement that
$(\ce, A)$ be generic in all our finiteness results however that the diagonals
should
be thrown away meaning that $P_{i}\neq P_{j}$ and $\gtb_{i}\neq
\gtb_{i}\neq\gtb_{j}$ for all $i\neq j$.

One may want to write down differential equations satisfied by horizontal
sections of this bundle. We will show in \ref{diffeqsatbycorrfuntttt} that
horizontal sections satisfy a system of $2m$ differential equations of
which
$m$equations are Knizhnik-Zamolodchikov equations and the other $m$ are
obtained from singular vectors of the Verma module projecting onto $M^{A}$.

Everything said here holds true for an admissible representation. It is easy to
see
that the bundle associated with an admissible representation is a quotient of
the just constructed bundle for the corresponding generalized Weyl module.

\subsection{Calculation of the dimensions of coinvariants. Fusion algebra}
Let $\ce\rightarrow\cp$ be the rank 2 trivial vector bundle and $M^{A}$ be a
$\hgtg^{A}-$module.
Here we will calculate the dimension of the space
$(M^{A})_{\gtg(\cp,A)},\;\sharp A=3$, in the
following 2 cases:
(i) the level $k$ is not
 rational and $M^{A}$ is a generalized Weyl module; (ii) $k+2=p/q$, $p\mbox{
and }q$
 being positive integers,
and $M^{A}$ is an admissible representation. Without loss of generality we can:

fix a coordinate $z$ on $\cp$; assume that
$A=\{(0,\gtb_{0}),\;(1,\gtb_{1}),\; (\infty,\gtb_{\infty})\}$, where
$\gtb_{0}=\nc e\oplus\nc h,\;\gtb_{\infty}=\nc f\oplus\nc h$ and
$\gtb_{1}=\nc(e-h-f)\oplus
\nc(h+2f)$.

(In fact, for any $\gtb_{0}\neq\gtb_{\infty}$  we can always choose a basis of
$\gtg$ so that
$\gtb_{0},\gtb_{\infty}$ are as above. As to $\gtb_{1}$, there really is some
freedom but it is
easy to see that all the calculations below are independent of the choice. We
have set
$\gtb_{1}=(\exp{f})\gtb_{0}(\exp{-f})$.)

\subsubsection{ The generic level case}
\label{Thegenericlevelcase}
So by \ref{defofocat}  we are given three irreducible generalized Weyl
 modules $V^{0}_{\lambda_{0},k},
V^{1}_{\lambda_{1},k},V^{\infty}_{\lambda_{\infty},k}$.
Recall, see \ref{Case_(i)},
that generalized Weyl modules are parametrized by symbols
$(V_{m}^{\epsilon},V_{n})$, where
$m,n$ are nonnegative integers,
$\epsilon\in\nz/2\nz$ and $V_{m}$ is an $m+1$-dimensional $\gtg-$module.
Therefore we
can and will assume that we have
\[(V_{m_{i}}^{\epsilon_{i}},V_{n_{i}}),\;i=0,1,\infty .\]

It is convenient to interpret the result of calculation of $\mbox{dim }
(\otimes_{i}(V_{m_{i}}^{\epsilon_{i}},V_{n_{i}}))_{\gtg(\cp,A)}$ in terms of
the
 {\em fusion algebra}.
The latter is defined as follows.   Suppose
 that for any pair of generalized Weyl modules, say
$(V_{r_{i}}^{\alpha_{i}},V_{s_{i}}),\;i=0,1$, there is only finite number of
$(V_{r_{\infty}}^{\alpha_{\infty}},V_{s_{\infty}})$ such that
\[\mbox{dim }
(\otimes_{i=0,1,\infty}(V_{r_{i}}^{\alpha_{i}},V_{s_{i}}))_{\gtg(\cp,A)}\neq
0.\]
Now view the symbols $(V_{m}^{\epsilon},V_{n})$ as generators of a free abelian
group. Then
there naturally arises an algebra (over $\nz$) with the operation of
  multiplication $\circ$ defined by
\[(V^{\alpha_{0}}_{r_{0}},V_{s_{0}})\circ (V^{\alpha_{1}}_{r_{1}},V_{s_{1}})=
\sum_{(r_{\infty},s_{\infty},\alpha_{\infty})}\mbox{dim }\{
(\otimes_{i=0,1,\infty}(V_{m_{i}}^{\epsilon_{i}},V_{n_{i}}))_{\gtg(\cp,A)}\}
(V^{\alpha_{\infty}}_{r_{\infty}},
V_{s_{\infty}}).\]
The algebra defined in this way is called {\em fusion
algebra}. Of course structure constants of the fusion
algebra determine the dimensions of the spaces of coinvariants.

One last piece  of notation:
 in the following
theorem  we  formally set $(X\oplus Y,Z)=(X,Z)+(Y,Z)$
and $(X,Y\oplus Z)=(X,Y)+(X,Z)$. Recall also that in the category of
$\gtg-$modules
one has
\[V_{r}\otimes V_{s}\approx V_{r+s}\oplus V_{r+s-2}\oplus\cdots\oplus
V_{|r-s|}.\]

\begin{theorem}
\label{fusalggencase}
(i) For any triple of generalized Weyl modules
the space $(V_{m_{i}}^{\epsilon_{i}},V_{n_{i}}))_{\gtg(\cp,A)}$ is finite
dimensional.

(ii) The fusion algebra is well-defined, multiplication being given by the
following formula
\begin{eqnarray}
(V_{r_{1}}^{\alpha},V_{s_{1}})\circ(V_{r_{2}}^{\beta},V_{s_{2}})&=&\nonumber\\
(V_{r_{1}+r_{2}}^{\alpha+\beta},V_{s_{1}}\otimes V_{s_{2}})&+&
(V_{r_{1}+r_{2}-1}^{\alpha+\beta+1},V_{s_{1}}\otimes V_{s_{2}})+
(V_{r_{1}+r_{2}-2}^{\alpha+\beta},V_{s_{1}}\otimes V_{s_{2}})+\cdots
+\nonumber\\
(V_{|r_{1}-r_{2}|}^{\alpha+\beta},V_{s_{1}}\otimes V_{s_{2}})&.&\nonumber
\end{eqnarray}
\end{theorem}

\subsubsection{ }
{\bf Proof of Theorem\ref{fusalggencase}.}

Throughout the proof $A$ will stand for $\{(0,\gtb_{0}),(1,\gtb_{1})\}$,
$A_{1}$ -- for
$\{(\infty,\gtb_{\infty})\}$. Along with the algebras
$\gtg(\cp,A),\; \gtg(\cp,A,A_{1})$ (see \ref{Finiteness_of_coinvariants})
introduce the
algebra $\bar{\gtg}(\cp,A,A_{1})\subset \gtg(\cp,A)$ consisting of all
functions
taking values in $\gtb_{\infty}$ at the point $\infty$.

 Of course
$\gtg(\cp,A,A_{1})\subset\bar{\gtg}(\cp,A,A_{1})$ is an ideal and
$\mbox{dim }\bar{\gtg}(\cp,A,A_{1})/\gtg(\cp,A,A_{1})=1$. Define
$\bar{h}_{\infty}$
to be a basis element of $\mbox{dim
}\bar{\gtg}(\cp,A,A_{1})/\gtg(\cp,A,A_{1})$.
It is a standard (and simple) fact of Lie algebra cohomology theory that
$\bar{h}_{\infty}$ acts on $(M^{A})_{\gtg(\cp,A,A_{1})}$.

\begin{lemma}
\label{cohom_for_Weyl}

 Let $M^{A}$ be a generalized Weyl module.  The element $\bar{h}_{\infty}$
 has a simple spectrum as an operator acting on
$(M^{A})_{\gtg(\cp,A,A_{1})} $. Further, if
$M^{A}=(V_{r_{1}}^{\alpha},V_{s_{1}})\otimes (V_{r_{2}}^{\beta},V_{s_{2}})$
then the
set of eigenvalues
of $\bar{h}_{\infty}$ is the set of the highest weights of the modules
appearing in
the right-hand side of Theorem\ref{fusalggencase}(ii).
\end{lemma}

Proof of this lemma is essentially the same as that of Theorem 4.4 in
\cite{feig_mal1} and
mostly consists of solving a system of 2 equations related to 2 singular
vectors -- one
in $(V_{r_{1}}^{\alpha},V_{s_{1}})$, another in
$(V_{r_{2}}^{\beta},V_{s_{2}})$.
We will discuss it in \ref{prooflemma49}. Derivation
of Theorem \ref{fusalggencase} from Lemma \ref{cohom_for_Weyl} is again very
similar
to that of Theorem 3.2 from Theorem 4.4 in {\em loc. cit}
and uses Verma modules as follows.

\begin{lemma}
\label{Berma-Weyl-gener}

(i) Let $(M^{A})^{\mu}_{\gtg(\cp,A,A_{1})}\in (M^{A})_{\gtg(\cp,A,A_{1})} $
 be the eigenspace related to the eigenvalue
$\mu$ of $\bar{h}_{\infty}$. Then $(M^{A})^{\mu}_{\gtg(\cp,A,A_{1})}\approx
(M^{A}\otimes M^{\infty,\gtb_{\infty}}_{\mu,k})_{\gtg(\cp,A\cup A_{1})}$.

(ii) Projection of a Verma module $M^{\infty,\gtb_{\infty}}_{\mu,k}$ onto a
generalized Weyl module $ W$ induces an isomorphism
of the coinvariants
\[(M^{A}\otimes M^{\infty,\gtb_{\infty}}_{\mu,k})_{\gtg(\cp,A\cup
A_{1})}\approx
(M^{A}\otimes W)_{\gtg(\cp,A\cup A_{1})}.\]
\end{lemma}

{\bf Proof of Lemma\ref{Berma-Weyl-gener}}

(i) A Verma module sitting at a point is induced from the 1-dimensional
representation of the
 algebra
of functions on the formal disk whose value at the point belong to the
corresponding Borel
subalgebra. Therefore (i) follows from Frobenius duality.

(ii)  Consider the resolution of $W$ by Verma modules (see \ref {Case_(i)},
formula
(\ref{exactseqforweyl}) ):
\[0\rightarrow M\rightarrow M^{\infty,\gtb_{\infty}}_{\mu,k}\rightarrow
W\rightarrow 0\]
and tensor it with $M^{A}$. There arises the long exact sequence of homology
groups of which
we consider the following part:

\[(M^{A}\otimes M)_{\gtg(\cp,A\cup A_{1})}\rightarrow
(M^{A}\otimes M^{\infty,\gtb_{\infty}}_{\mu,k})_{\gtg(\cp,A\cup
A_{1})}\rightarrow
(M^{A}\otimes W)_{\gtg(\cp,A\cup A_{1})}\rightarrow 0.\]

Since $M^{\infty,\gtb_{\infty}}_{\mu,k}$ projects onto a Weyl module, the Verma
module $M$
does not, see \ref{Case_(i)}. Lemma\ref{cohom_for_Weyl} and  now give that
$(M^{A}\otimes M)_{\gtg(\cp,A\cup A_{1})}=\{0\}.$ $\Box$

\bigskip

To complete the proof of Therorem\ref{fusalggencase} observe that
Lemma \ref{cohom_for_Weyl}  and Lemma \ref{Berma-Weyl-gener} together is a
reformulation of
Therorem\ref{fusalggencase}. $\Box$
\bigskip

\begin{corollary} Let $\sharp A=1$ and let $A_{1}$ and $A_{2}$
be as in \ref{Finiteness_of_coinvariants}. The following conditions are
equivalent
\label{oncompseries-gener}

(i) $M^{A}$ is a direct sum of generalized Weyl module;

(ii)$SSM^{A}=\Omega^{A}_{reg}\cup\sigma\Omega^{A}_{reg}$;

(iii) For any Verma module $W^{A_{1}}$
$\mbox{dim}(M^{A}\otimes W^{A_{1}})_{\gtg(\cp,A\cup A_{1},A_{2})}<\infty$.
\end{corollary}

\subsubsection{ }
\label{prooflemma49}
Here we sketch the proof of Lemma \ref{cohom_for_Weyl}. First of all replace
$M^{A}$ with
the corresponding Verma module -- $\bar{M}^{A}$. Then
pass from the space
$(\bar{M}^{A})_{\gtg(\cp,A,A_{1})}$ to its dual, that is to the space
of $\gtg(\cp,A,A_{1})-$invariant functionals on $\bar{M}^{A}$.
Choose $h\otimes (1-z^{-1})$ to be a representative of $\bar{h}_{\infty}$. Let
$\Psi$ be the
 eigenvector of $h\otimes (1-z^{-1})$. By definition $\Psi$ is a linear
functional on
$M_{\lambda_{0},k}^{0,\gtb_{0}}\otimes M_{\lambda_{1},k}^{1,\gtb_{1}}$. It is
an excersise
on Frobenius duality to show that such a functional exists and unique.

 Define $F$ to be the
following linear
functional on $M_{\lambda_{0},k}^{0,\gtb_{0}}$: $F(w)=\Psi(w\otimes
v_{\lambda_{1}})$, where,
as usual, $v_{\lambda_{1}})$ is the vacuum vector of
$M_{\lambda_{1},k}^{1,\gtb_{1}}$. As
$M_{\lambda_{0},k}^{0,\gtb_{0}}$ is $\nz_{+}\times\nz_{+}-$graded (see
\ref{defofalgebrandmodusu}),
we denote by $F_{ij}$ the restriction of $F$ to the $(i,j)-$component. Direct
calculations show that
with respect to the natural action of $\hgtg$ on
$M_{\lambda_{0},k}^{0,\gtb_{0}}$:

\begin{eqnarray}
\label{strofloopmodact}
& &\oplus_{i,j\in\nz}\nc F_{ij}\approx \cf_{\alpha\beta}^{\nc^{\ast}}\\
& &\mbox{where }\alpha=\frac{\lambda_{\infty}-\lambda_{1}-\lambda_{0}-2}{2},\;
\beta=\lambda_{1}
\end{eqnarray}

The functional $F$ factors through the projection
$M_{\lambda_{0},k}^{0,\gtb_{0}}\rightarrow V_{\lambda_{0},k}^{0,\gtb_{0}}$ if
and only if
it vanishes on the singular vector of $M_{\lambda_{0},k}^{0,\gtb_{0}}$. In
other words, if this
singular vector, say $S$, has degree $(i,j)$ then the following equation holds
\[ S F_{ij}=0.\]
The latter equation can be  written down and solved explicitly using formulas
(\ref{actonloops1} or
\ref{actonloops2}). Similar arguments go through for the module
$M_{\lambda_{1},k}^{1,\gtb_{1}}$
giving another equation, say
\[ S' F_{i'j'}=0.\]

Simultaneous solutions to these 2 equations give the desired result.
By the way, as (\ref{actonloops1},
\ref{actonloops2}) show, each of the expressions $ S F_{ij},\; S' F_{i'j'}$
splits in a product of linear factors; therefore geometrically the solution is
a collection of intersection points of 2 families of lines in the plane. $\Box$

\subsubsection {The rational level case}
\label{The_rational_leve_case}
Suppose $k+2=p/q$, $p\mbox{ and }q$
 being positive integers.
Now instead of 3 generalized Weyl modules sitting at 3 points in $\cp$ we are
given
3 admissible representations sitting at 3 points on $\cp$.
Recall, see \ref{Case_(ii)}, that  admissible representations
are parametrized by symbols $(V_{m}^{\epsilon},V_{n}),\;0\leq m\leq q-1,0\leq
n\leq p-2$ modulo the relation
$(V_{m}^{\epsilon},V_{n})=(V_{q-1-m}^{\epsilon+1},V_{p-2-n})$. Denote
by $(V_{m}^{\epsilon},V_{n})^{\sim }$  an equivalence class
of $(V_{m}^{\epsilon},V_{n})$. We assume that $(V_{m}^{\epsilon},V_{n})^{\sim}$
satifies
the same bilinear condition  $(V_{m}^{\epsilon},V_{n})$ in
Theorem\ref{fusalggencase} does.

The definition of the fusion (Verlinde) algebra in this case repeats
word for word that in \ref{Thegenericlevelcase}.

Recall finally that  Kazhdan-Lusztig fusion functor \cite{kazh_luszt} gives
\[V_{r}\dot{\otimes}_{k}V_{s}= V_{|m-n|}
\oplus V_{|m-n|+2}\cdots \oplus V_{\mbox{min}\{2k-r-s,r+s\}}.\]
The following theorem was proved in \cite{feig_mal1}
in an equivalent but much less illuminating form.

\begin{theorem}
\label{fusalgratcase}
(i) For any triple of admissible representations
the space $(V_{m_{i}}^{\epsilon_{i}},V_{n_{i}}))_{\gtg(\cp,A)}$ is finite
dimensional.

(ii) The fusion algebra is well-defined, multiplication being given by the
following formula
\begin{eqnarray}
(V_{r_{1}}^{\alpha},V_{s_{1}})^{\sim }\circ(V_{r_{2}}^{\beta},V_{s_{2}})^{\sim
}&=&\nonumber\\
(V_{|r_{1}-r_{2}|}^{\alpha+\beta},V_{s_{1}}\dot{\otimes}_{p-2} V_{s_{2}})^{\sim
}&+&
(V_{|r_{1}-r_{2}|+1}^{\alpha+\beta},V_{s_{1}}\dot{\otimes}_{p-2}
V_{s_{2}})^{\sim }+
(V_{|r_{1}-r_{2}|+2}^{\alpha+\beta},V_{s_{1}}\dot{\otimes}_{p-2}
V_{s_{2}})^{\sim }+\cdots +
\nonumber\\
(V_{N}^{\alpha+\beta},V_{s_{1}}\dot{\otimes}_{p-2} V_{s_{2}})^{\sim }&
&\nonumber,
\end{eqnarray}
where $N=\mbox{min}\{2q-2-r-s,r+s\}$.
\end{theorem}
 It is an easy exercise to derive this theorem from Theorem
\ref{fusalggencase}. For future
purposes, however, we now sketch its original proof.  Set
$A=\{(\infty,\gtb_{\infty}\}$,
$A_{2}=\{(0,\gtb_{0}),(1,\gtb_{1})\}$. In addition to the algebras
$\gtg(\cp,A),\gtg(\cp,A,A_{2})$
as in \ref{Finiteness_of_coinvariants}, we introduce an algebra
$\bar{\gtg}(\cp,A,A_{2})\subset\gtg(\cp,A)$. The latter consists of all
functions whose values
at the points $0$ ($1$ resp.) belong to $\gtb_{0}$ ($\gtb_{1}$ resp.).
Obviously
$\gtg(\cp,A,A_{2})\subset\bar{\gtg}(\cp,A,A_{2})$ is an ideal and the quotient
algebra
$\bar{\gtg}(\cp,A,A_{2})/\gtg(\cp,A,A_{2})$ is commutative and 2-dimensional.
This algebra naturally operates on the space
$(M^{A})_{\gtg(\cp,A,A_{2})}$. Let $\bar{h}_{0},
\bar{h}_{1}$ be a basis of $\bar{\gtg}(\cp,A,A_{2})/\gtg(\cp,A,A_{2})$.

\begin{lemma}(\mbox{\cite{feig_mal1}})
\label{cohom_for_admiss}

 (i) $\mbox{dim}(M^{A})_{\gtg(\cp,A,A_{2})}<\infty$ if and only if $M^{A}$ is
  an admissible representation.

(ii) Let $M^{A}$ be an admissible representation. The elements
$\bar{h}_{0},
\bar{h}_{1}$
 have simple spectra
as  operators acting on
$(M^{A})_{\gtg(\cp,A,A_{2})}$. Their eigenvalues recover the structure
constants of the
fusion algebra.
\end{lemma}

``Inserting'' Verma modules and using BGG resolution one derives Theorem
\ref{fusalgratcase}
from Lemma\ref{cohom_for_admiss} in  a way similar to that we used in
\ref{Thegenericlevelcase}.

Another important corollary of Lemma\ref{cohom_for_admiss} is as follows.

\begin{corollary} Let $\sharp A=1$ The following conditions are equivalent
\label{oncompseries-ratin}
(i) $M^{A}$ is a sum of admissible representations;

(ii)$SSM^{A}=\Omega^{A}_{nilp}\cup\sigma\Omega^{A}_{nilp}$;

(iii)
$\mbox{dim}(M^{A})_{\gtg(\cp,A,A_{2})}<\infty$.
\end{corollary}

\subsubsection {Classical and quantum $\gtosp(1|2)$. Fusion algebra as a
Grothendieck ring}
\label{fusalgasgrotring}

{\bf A. } $\gtosp(1|2)$ is a rank 1 superalgebra -- one of the superanalogues
of $\gtsl_{2}$. It can be defined as an algebra on 2 odd generators,
$x_{+},x_{-}$, one even generator, $h$, and relations
\[[x_{+},x_{-}]=h,\; [h,x_{\pm}]=\pm x_{\pm}.\]

Even part of this algebra is $\gtsl_{2}$ and is generated by $x_{\pm}^{2}$;
odd part is $V_{1}$ as an $\gtsl_{2}-$ module, its basis is $x_{+},\; x_{-}$.

 From this it is easy to obtain the following classification of all simple
finite dimensional $\gtosp(1|2)-$modules. (It is even simpler to do this in the
way
modelling the $\gtsl_{2}-$case -- by starting with Verma modules and then
quotienting out a singular vector; for details see \cite{kul_resh}). Each
$\gtosp(1|2)-$module  $W$is a sum of an even and odd part
$W=^{even}W\oplus^{odd}W$;
each $^{\cdot}W$ is an $\gtsl_{2}-$module, i.e. direct sum of $V_{n}$'s. These
are generalities. But in reality each irreducible $\gtosp(1|2)-$module is of
one of the 2 following
types:
\[ V_{n}^{0} \mbox{ such that } ^{even}V_{n}^{0}=V_{n},\;
^{odd}V_{n}^{0}=V_{n-1};\]

\[ V_{n}^{1} \mbox{ such that } ^{even}V_{n}^{0}=V_{n-1},\;
^{odd}V_{n}^{0}=V_{n}.\]

The fact that the dimensions of the even and odd parts are different by 1 is a
consequence  of the fact  that odd part of the algebra is $V_{1}$.

We see that each irreducible $\gtosp(1|2)-$module is odd-dimesional; further
$V_{n}^{0}$ and $V_{n}^{1}$ are isomorphic as modules and obtained from each
other by the change of parity. This is the category of finite dimensional
representations of $\gtosp(1|2)$; denote it $Rep(\gtosp(1|2))$. As in the
$\gtsl_{2}-$case, one proves that $Rep(\gtosp(1|2))$ is semisimple.

The universal enveloping algebra $U\gtosp(1|2)$ is in fact a Hopf algebra, for
example the comultiplication is given by the standard formula $g\mapsto
g\otimes 1+1\otimes g,\; g\in\gtosp(1|2)$. This makes $Rep(\gtosp(1|2))$ a
tensor category: $Rep(\gtosp(1|2))\times Rep(\gtosp(1|2))\rightarrow
Rep(\gtosp(1|2))$, $A,B\mapsto A\otimes B$, where the $\gtosp(1|2)-$module
structure on $A\otimes B$ is determined through the comultiplication (and the
rule of sign!). Decomposing the tensor product of 2 irreducible modules one
gets the Grothendieck ring of $Rep(\gtosp(1|2))$.

\begin{lemma}
\label{decomposptensprod}
\[ V_{r_{1}}^{\alpha}\otimes V_{r_{2}}^{\beta}=
V_{r_{1}+r_{2}}^{\alpha+\beta}+
V_{r_{1}+r_{2}-1}^{\alpha+\beta+1}+
V_{r_{1}+r_{2}-2}^{\alpha+\beta}+\cdots +
V_{|r_{1}-r_{2}|}^{\alpha+\beta}.\]
\end{lemma}
{\bf Proof.} Direct calculations show that $ V_{r_{1}}^{\alpha}\otimes
V_{r_{2}}^{\beta}$ contains one and only one singular (annihilated by $x_{+}$)
vector of each weight from $|r_{1}-r_{2}|$ to $r_{1}+r_{2}$ and that the
submodules generated by these vectors are irreducible. Proof is completed by
counting dimensions. $\Box$

\bigskip

 Theorems \ref{fusalggencase} and \ref{fusalgratcase} provide us with 2
commutative algebras. Here
we interprete these algebras as Grothendieck rings of certain categories. Start
with the algebra
of Theorem \ref{fusalggencase} and denote it $\ca^{gen}$. Obviously
$\ca^{gen}=\ca_{0}\otimes\ca$, where $\ca$ is the Grothendieck ring of the
category of
finite-dimensional representations of $\gtg$ (its multiplication law is defined
by the formula
preceding Theorem \ref{fusalggencase}) and $\ca_{0}$ is the algebra with basis
$V_{i}^{\alpha},\;i\geq 0,\alpha\in\nz/2\nz$, multiplication being given by
\begin{equation}
\label{reprospalg}
V_{r_{1}}^{\alpha}\circ V_{r_{2}}^{\beta}=
V_{r_{1}+r_{2}}^{\alpha+\beta}+
V_{r_{1}+r_{2}-1}^{\alpha+\beta+1}+
V_{r_{1}+r_{2}-2}^{\alpha+\beta}+\cdots +
V_{|r_{1}-r_{2}|}^{\alpha+\beta}.
\end{equation}

Comparing (\ref{reprospalg}) with Lemma \ref{decomposptensprod} we get the
following.

\begin{proposition}
\label{interptofgenlev}
$\ca_{0}$ is the Grothendieck ring of the category of finite-dimensional
representations
of the superalgebra $ \gtosp(1|2)$.
\end{proposition}

Appearance of $\gtosp(1|2)$ here, although artificial as it may seem to be, has
deep reasons  behind it. To see this we will analyze the rational level case
using quantized enveloping algebras.

\begin{remark}
It follows from Lemma \ref{decomposptensprod} that the functor
$Rep(\gtosp(1|2))\rightarrow Rep(\gtsl_{2})$, $V_{m}^{\alpha}\mapsto
V_{m}\oplus V_{m-1}$ induces  an epimorhpism of the Grothendieck rings.
\label{renonmorh[pring}
\end{remark}

 \bigskip

\bigskip

{\bf B.} Both $U\gtsl_{2}$ and $U\gtosp(1|2)$ admit quantization,
$U_{t}\gtsl_{2}$ and $U_{t}\gtosp(1|2)$ resp.. Let us remind the relevant
formulas.
The Drinfeld-Jimbo (see \cite{drinf,jimbo}) algebra $U_{t}\gtsl_{2},\; t\in\nc$
is defined to be an associative
algebra on generators $E,F,K^{\pm 1}$ and relations

\[EF-FE=\frac{K-K^{-1}}{t-t^{-1}},\; KEK^{-1}=t^{2}E,\; KFK^{-1}=t^{-2}F.\]

$U_{t}\gtosp(1|2)$ is similarly defined \cite{kul_resh} as an associative
algebra on generators $X_{+}, X_{-}, K^{\pm 1}$ and relations

\[X_{+}X_{-}+X_{-}X_{+}=\frac{K-K^{-1}}{t-t^{-1}},\; KX_{\pm}K^{-1}=t^{\pm
1}X_{\pm}.\]

The representation theory of $\gtsl_{2}$ and $\gtosp(1|2)$ ``deforms to'' the
representation theory of $U_{t}\gtsl_{2}$ and  $U_{t}\gtosp(1|2)$ resp. We will
 continue denoting by $V_{m}$ the $m+1-$dimensional module over
$U_{t}\gtsl_{2}$, and by $V^{0}_{m},\; V^{1}_{m}$ the 2 $(2m+1)-$ dimensional
modules over $U_{t}\gtosp(1|2)$. For generic $t$ these modules are irreducible,
the categories of finite dimensional representations, $Rep(U_{t}\gtsl_{2})$ and
$Rep(U_{t}\gtosp(1|2))$, generated by these modules are semisimple.

The deformations $U_{t}\gtsl_{2}$ and $U_{t}\gtosp(1|2)$ are especially
remarkable in that they afford simultaneous deformation of the Hopf algebra
structure. We get 2 tensor categories $Rep(U_{t}\gtsl_{2})$ and
$Rep(U_{t}\gtosp(1|2))$. What has been said implies that the Grothendieck rings
of $Rep(U_{t}\gtsl_{2})$ and
$Rep(U_{t}\gtosp(1|2))$ are isomorphic to the Grothendieck rings of the
corresponding classical objects.

If however $t$ is a root of unity, things change dramatically. Suppose for
simplicity that $t$ is a primitive $l$-th root of unity, $l$ being odd. Then

\begin{equation}
\mbox{(i) $V_{m}$ is irreducible if and only if $m< l$;}
\label{irratroots1}
\end{equation}
\begin{equation}
\label{irratroots2}
\mbox{(ii) $V_{m}^{\epsilon}$ is irreducible if and only if $m< l$.}
\end{equation}

(Both statements are proved by direct computations.)

What is even more important is that the categories $Rep(U_{t}\gtsl_{2})$ and
$Rep(U_{t}\gtosp(1|2))$ are no longer semisimple. For example, tensor product
of 2 irreducible representations is not semisimple. Things, however, are still
very
much under control.

\begin{lemma}
\label{decompprodatroots}
Let $t$ be a primitive $l$-th root of unity, $l$ being odd, $m,n<l$. Then

\[\mbox{ (i) } V_{m}\otimes V_{n}= V_{|m-n |}
\oplus V_{|m-n|+2}\cdots \oplus V_{\mbox{min}\{2(l-1)-m-n,m+n\}}\oplus W,\]
where $W$ is not semisimple.
\[\mbox{ (ii) } V_{m}^{\alpha}\otimes V_{n}^{\beta}= V_{|m-n |}^{\alpha+\beta}
\oplus V_{|m-n|+1}^{\alpha+\beta+1}\cdots \oplus
V_{\mbox{min}\{2(l-1)-m-n,m+n\}}^{\alpha+\beta}\oplus W,\]
where $W$ is not semisimple.

\end{lemma}

{\bf Sketch of Proof.} (i) is well-known, see \cite{resh_tur}. We will however
review both cases as at our level of brevity there will no difference between
them. First, direct calculations as in the proof of Lemma
\ref{decomposptensprod} show that  regardless of $t$ at each weight space there
can always be only one singular vector. Now decomposition of Lemma
\ref{decomposptensprod}, statements (\ref{irratroots1}, \ref{irratroots2})
and this uniqueness result show that the submodules $V_{l-1+i}$ and $V_{l-1-i}$
(or $V_{l-1+i}^{\epsilon}$ and $V_{l-1-i}^{\epsilon}$), $i\leq m+n-l+1$ are
non-trivially tangled. Other $V_{j}$ coming from generic $t$ are still
irreducible and appear as direct summands. $\Box$

\bigskip

{\bf Definition.}

(i)  Define $Rep(U_{t}\gtsl_{2})^{(l)}$ and
$Rep(U_{t}\gtosp(1|2))^{(l)}$ to be subcategories of $Rep(U_{t}\gtsl_{2})$ and
$Rep(U_{t}\gtosp(1|2))$resp. consisting of direct sums of irreducible modules
$V_{m}$ (or $V_{m}^{\alpha}$ resp.), $m< l$.

(ii) Define functors
\[Rep(U_{t}\gtsl_{2})^{(l)}\times Rep(U_{t}\gtsl_{2})^{(l)}\rightarrow
Rep(U_{t}\gtsl_{2})^{(l)},\; A,B\mapsto A\dot{\otimes}B,\]

\[Rep(U_{t}\gtosp(1|2))^{(l)}\times Rep(U_{t}\gtosp(1|2))^{(l)}\rightarrow
Rep(U_{t}\gtosp(1|2))^{(l)},\; A,B\mapsto A\dot{\otimes}B,\]
by taking the usual tensor product and then throwing away $W$ in the right hand
side of formulas in Lemma \ref{decompprodatroots}. $\Box$

\bigskip

We get tensor categories $Rep(U_{t}\gtsl_{2})^{(l)}$ and
$Rep(U_{t}\gtosp(1|2))^{(l)}$.

\bigskip

\bigskip

{\bf C.} It is easy now to interpret the  fusion algebra at the rational level
in terms of the Grothendieck rings of $Rep(U_{t}\gtsl_{2})^{(l)}$ and
$Rep(U_{t}\gtosp(1|2))^{(l)}$.   In view of Lemma \ref{decompprodatroots} and
Definition
above, Theorem \ref{fusalgratcase} reads as follows.

\begin{proposition}
\label{fusalgeasgrotratlev}
Fusion algebra at the level $k+2=p/q$ is a quotient of the tensor product
of the Grothendieck rings of the categories  $Rep(U_{t_{1}}\gtsl_{2})^{(p-1)}$
and
$Rep(U_{t_{2}}\gtosp(1|2))^{(q)}$. Further, the fusion algebra always contains
the Grothendieck ring of $Rep(U_{t_{2}}\gtosp(1|2))^{(q)}$ via the classes
of symbols $V_{m}^{\alpha}, V_{0}$.
\end{proposition}

\subsubsection{Kac-Moody vs. Virasoro}
\label{kacmoodyvsvir}
Virasoro algebra, $Vir$, is defined to be a vector space with basis
$\{ L_{i},z, \; i\in\nz\}$ and bracket
\[[L_{i},L_{j}]=(i-j)L_{i+j}+\delta_{i,-j}\frac{i^{3}-i}{12}z.\]
Representation theory of Virasoro algebra is to a great extent parallel
to that of $\hgtg$. We will confine to essentials, making reference to
\cite{fei_fuchs}.

One defines the Verma module $M_{h,c}$, where $(h,c)$ is a highest weight, i.e.
eigenvalues of $L_{0}\;, z$ resp. determined by the vacuum vector; $c$ is
sometimes referred to as level. A Verma module
is reducible if and only if it contains a singular vector.  $M_{h,c}$
generically has no singular vectors.By the Kac determinant formula, there is a
family of hyperbolas labelled by pairs of positive integers $m,n$ in the plane
with coordinates $(h,c)$ such that if $M_{h,c}$ contains a singular vector,
then $(h,c)$ belongs to one of  these hyperbolas; generically along hyperbolas
the singular vector is unique.   Denote the singular vector arising in
$M_{h,c}$ as
$(h,c)$ gets on hyperbola with the label $m,n$ by $S_{mn}$. There arises
the $Vir$-analogue of the generalized Weyl module $M_{h,c}/<S_{mn}=0>$.
Attach to        $M_{h,c}/<S_{mn}=0>$ the symbol $(V_{n-1},V_{m-1})$. Further,
for $c$ fixed there arises a one-to-one correspondence between the
$Vir$-analogues of generalized Weyl modules and symbols $(V_{n-1},V_{m-1})$.
This has all been in precise analogy with \ref{Case_(i)}.

It has hardly been written anywhere, but is nevertheless known that the $Vir-$
analogue
of the fusion algebra from \ref{Thegenericlevelcase}, i.e. at a generic level,
is as follows:
\begin{equation}
\label{fusalgforvir}
(V_{n_{1}},V_{m_{1}})\circ (V_{n_{2}},V_{m_{2}})=(V_{n_{1}}\otimes
V_{n_{2}},V_{m_{1}}\otimes V_{m_{2}}).
\end{equation}
(The interested reader can prove this result using methods of
\cite{feig_fuchs_2}; our treatment of the $\hgtg$-fusion algebra in
\ref{Thegenericlevelcase} is also a direct analogue of these.)

There is a functor sending $\hgtg-$modules to $Vir-$modules -- quantum
Drinfeld-Sokolov reduction. One of the prerequisites for it is a choice of
a nilpotent subalgebra of $\gtsl_{2}$. The two obvious possibilities are
$\nc e$ and $\nc f$. Denote the corresponding functors $\phi_{e}$ and
$\phi_{f}$. It can be extracted from \cite{feig_fr} that  both functors send
generalized Weyl modules to generalized Weyl modules. In our terminology one
gets

\[\phi_{e}:\begin{array}{lll}
(V_{m}^{0},V_{n})&\mapsto& (V_{m},V_{n})\\
(V_{m}^{1},V_{n})&\mapsto& (V_{m-1},V_{n}),
\end{array}\]

\[\phi_{f}:\begin{array}{lll}
(V_{m}^{0},V_{n})&\mapsto& (V_{m-1},V_{n})\\
(V_{m}^{1},V_{n})&\mapsto& (V_{m},V_{n}),
\end{array}\]
where the symbol $V_{-1}$, if arises, is understood as zero.

The $Vir-$analogue of admissible representations is the celebrated minimal
representations. The latter can be defined as quotients of generalized Weyl
modules
by repeating word for word definition of  admissible representations from
\ref{Case_(ii)}. It is known that minimal representations
arise only when
 \[c=c_{pq}=1-\frac{6(p-q)^{2}}{pq},\]
where $p,q$ are relatively prime positive integers.
 There are again 2 generalized Weyl modules projecting on
a given minimal representation. Therefore minimal representations are labelled
by equivalence classes of symbols $(V_{m}, V_{n})$. It can be shown that the
equivalence relation is as follows: $(V_{m},V_{n})\approx (V_{q-2-m},
V_{p-2-n})$ for $c=c_{pq}$. From this and (\ref{fusalgforvir}) one can easily
calculate the fusion algebra. We will not write down the relevant formulas here
and confine to mentionaing that the algebra is related to the product of
Grothendieck rings of 2 quantum $U_{t}(\gtsl_{2})$ at appropriate roots of
unity in much the same way as the fusion algebra for $\hgtg$ is related to the
product
of Grothendieck rings of  $U_{t}(\gtosp(1|2))$ and $U_{t}(\gtsl_{2})$. Recall
also that the $Vir$-fusion algebra  was calculated in \cite{bpz};
mathematically acceptable exposition can be found in \cite{feig_fuchs_2}.

Another property of the Drinfeld-Sokolov reduction is that both  $\phi_{e}$ and
$\phi_{f}$ send admissible representations at the level $k=2-p/q$ of $\hgtg$ to
minimal representations
of $Vir$ at the level $c_{pq}$, see \cite{fkw} .

\begin{proposition}
\label{homosl2vir}
 The functor $\phi_{e}\oplus\phi_{f}$ determines an epimorphism of the
$\hgtg-$fusion algebra onto the $Vir-$fusion algebra at both generic and
rational levels.
\end{proposition}

{\bf Proof.} The generic level case follows from Remark \ref{renonmorh[pring}
and formula (\ref{fusalgforvir}) above. In the rational level case, the
statement follows from the fact that both, $\hgtg-$ and $Vir$-, fusion algebras
are obtained from their generic level counterparts by imposing the equivalence
relations and the 2 equivalence relations agree with each other.  $\Box$

\subsection{Fusion functor.}
This part is an announcement, proofs will appear elsewhere

Suppose we have a trivial vector bundle $\ce\rightarrow\cp$,
$A=\{(P_{1},\gtb_{1}),(P_{2},\gtb_{2})\}$, $B=\{(P_{3},\gtb_{3})\}$, so that
$(\ce,A\bigsqcup B)$ is generic. There is a construction which to a
$\hgtg^{A}-$module associates a $\hgtg^{B}-$module. This construction is a
natural adjustment of the Kazhdan-Lusztig tensoring \cite{kazh_luszt} to our
needs.

Denote by $\gtg(\cp,A,B)$ the subalgebra of $\gtg(\cp,A)$ consisting of
functions taking values in $\gtn_{3}=[\gtb_{3},\gtb_{3}]$ -- just like we did
in \ref{finofcoinvsferesubs}. For a $\hgtg^{A}-$module $M^{A}$, denote by
$M^{A}_{N}$ the subspace of $(M^{A})^{\ast}$ annihilated by
$\gtg(\cp,A,B)^{N}$. Obviously  $M^{A}_{N}\subset M^{A}_{N+1}$, $N\geq 1$. Set

\[F^{A\rightarrow B}(M^{A})=\cup_{N\geq 1} M^{A}_{N}.\]

One can show that the vector space
 $F^{A\rightarrow B}(M^{A})$ affords in a natural way a structure of an
$\hgtg^{B}-$module at the same level; this is easy to show in the spirit of
\cite{kazh_luszt, beil_feig_maz}. Using our methods one can show that

(i) if $M^{A}$ is from the $\co-$category, or further a generalized Weyl
module, or further an admissible representation, then $F^{A\rightarrow
B}(M^{A})$ is also as a $\hgtg^{B}-$module;

(ii) the arising in this way Grothendieck rings coincide with those in Theorem
\ref{fusalggencase} or Theorem
\ref{fusalgratcase} if the level is generic or rational resp..

This generalizes the statement for the integrable representations, see
\cite{fink}.

{\bf Problem.} Describe the arising tensoring in the spirit of Kazhdan-Lusztig.

 \subsection{Quadratic degeneration}
\label{Quadratic_degeneration}

\subsubsection { }
\label{setupquadrdeg}
The setup here will the following version of  \ref{maingenerresults}:

(i) $\bar{\pi}:\;\cc_{S}\rightarrow S$ be a family of  curves over  a formal
disk  $S$, such that the fiber over the generic point of $S$ (``outside
origin'')
is a smooth projective curve, and over the origin, $O$, the fiber is a curve
$\cc_{O}$.
with exactly one quadratic singularity;

(ii) $\rho_{S}:\; \ce_{S}\rightarrow \cc_{S}$ is a rank $2$ vector bundle.

As in \ref{maingenerresults}, we complete these data to the localization data
with logarithmic singularities, say
$\tilde{\psi}$.
In the standard way, Theorem \ref{existofDmod} rewrites to give a $D-$module
over $S$ with logarithmic singularities at $O$; call it
$\Delta_{\tilde{\psi}}(M^{A})$. This is because $Spec(S)$ is $\nc[[t]]$ and
vector fields
vanishing at $q=0$ are exactly those which can be lifted to $\cc_{S}$.

Along with the family $\bar{\pi}:\;\cc_{S}\rightarrow S$ consider the family
$\bar{\pi}^{\vee}:\;\cc_{S}^{\vee}\rightarrow S$, obtained from
$\bar{\pi}:\;\cc_{S}\rightarrow S$ by replacing the singular fiber $\cc_{O}$
with
its normalization $\cc_{O}^{\vee}$ (i.e. be tearing $\cc_{O}$ apart at the
self-intersection point). There is a projection $\cc_{O}^{\vee}\rightarrow
\cc_{O}$
and the preimage of the self-intersection point $a\in C_{O}$ consists
of 2 points $a_{0}, a_{\infty}\in C_{O}^{\vee}$.

It is obvious that the datum $\ce\rightarrow \cc_{S}$
is equivalent to the data  ``$\rho^{\vee}_{S}: \ce_{S}^{\vee}\rightarrow
\cc_{S}^{\vee}$, equivalence $(\rho^{\vee}_{S})^{-1}(a_{0})\approx
(\rho^{\vee}_{S})^{-1}(a_{\infty})$''. The localization data with logarithmic
singularities $\tilde{\psi}$ rewrites to give a ``normalized'' localization
data
$\psi^{\vee}$.

In addition fix 2 different lines $l_{0},l_{\infty}$ in the fiber of $\ce_{S}$
over the point $a\in\cc_{O}$. This determines 2 Borel subalgebras,
$\gtb_{0},\gtb_{\infty}$ operating in the fiber over $a$.

After normalization these additional data determine the line $l_{0}$ and the
Borel subalgebra $\gtb_{0}$ operating in the fiber of $\ce_{S}^{\vee}$ over
$a_{0}$, as well as the line $l_{\infty}$ and the Borel subalgebra
$\gtb_{\infty}$ operating in the fiber over $a_{\infty}$. We also get a
distinguished Cartan subalgebra
$\gth=\gtb_{0}\cap\gtb_{\infty}$. Set
$A^{\vee}=A\bigsqcup\{(a_{0},\gtb_{0}),(a_{\infty},\gtb_{\infty})\}$.

Now with a $\hgtg^{A}-$module $M^{A}$ at the level $k$ and an admissible weight
  $\lambda\in\gth^{\ast}$ we associate the $\hgtg^{A^{\vee}}-$module
$M^{A}\otimes L^{P_{0},\gtb_{0}}_{\lambda,k}\otimes
L^{P_{\infty},\gtb_{\infty}}_{\lambda,k}$.  We get a $D-$module for the
``normalized''localization data:

 \[\oplus_{\lambda}\Delta_{\psi^{\vee}}(M^{A}\otimes
L^{P_{0},\gtb_{0}}_{\lambda,k}\otimes
L^{P_{\infty},\gtb_{\infty}}_{\lambda,k}).\]

\begin{proposition}
\label{statemonquadrgen}
Generically with respect to $l_{0},l_{\infty}$,
if $\Delta_{\tilde{\psi}}(M^{A})$ is smooth then
$\oplus_{\lambda}\Delta_{\psi^{\vee}}(M^{A}\otimes
L^{P_{0},\gtb_{0}}_{\lambda,k}\otimes
L^{P_{\infty},\gtb_{\infty}}_{\lambda,k})$ is also and there is an isomorphism
of $D-$modules
\[\Delta_{\tilde{\psi}}(M^{A})\approx
\oplus_{\lambda}\Delta_{\psi^{\vee}}(M^{A}\otimes
L^{P_{0},\gtb_{0}}_{\lambda,k}\otimes
L^{P_{\infty},\gtb_{\infty}}_{\lambda,k}).\]
\end{proposition}

\subsubsection { Proof}

(i) Begin with the genus zero case.Observe that the algebra of regular
functions
on the neighborhood of the point $a$ is
$\nc[t_{0},t_{\infty}][[t]]/<t_{0}t_{\infty}=t>$ where $t$ is a coordinate on
$S$; $\cc_{O}^{\vee}$
in this case is just a union of 2 spheres. Therefore the set $A$ splits in two:
$A'$ and $A''$, each of which has to do with one of the spheres.

Hence the algebra $\gtg(\bar{\pi}^{-1}(s),A)$ can be degenerated into the
following one
as $s$ ``approaches'' $O$:

 \[(\gtg(\cp,A',(P_{0},\gtb_{0}))+\gth)\oplus_{\gth}
(\gth+\gtg(\cp,A'',(P_{\infty},\gtb_{\infty})).\]

Meaning of the last expression is as follows: recall, see
\ref{Finiteness_of_coinvariants}, that $\gtg(\cp,A',(P_{0},\gtb_{0}))$ consists
of functions regular outside $\bar{A}$ and sending $P_{0}$ to $\gtn_{0}$;
$\gtg(\cp,A'',(P_{\infty},\gtb_{\infty}))$ is defined similarly with
$P_{0},\gtn_{0}$ replaced with $P_{\infty},\gtn_{\infty}$; further the algebra
$\gtg(\cp,A',(P_{0},\gtb_{0}))+\gth$ is the algebra of functions sending
$P_{0}$ to $\gtb_{0}$, the same is true for
$\gth+\gtg(\cp,A'',(P_{\infty},\gtb_{\infty}))$; finally ``$\oplus_{\gth}$''
means
direct product over $\gth$.

Therefore the coinvariants degenerate into the space
\[((M^{A'})_{\gtg(\cp,A',(P_{0},\gtb_{0}))}\otimes
(M^{A''})_{\gtg(\cp,A'',(P_{\infty},\gtb_{\infty}))})_{\gth},\]
where $\gth$ acts by means of the diagonal embedding; this makes sense as the
fibers are identified.

By Proposition \ref{finofcoinv_general}, the space
\[(M^{A'})_{\gtg(\cp,A',(P_{0},\gtb_{0}))}\otimes
(M^{A''})_{\gtg(\cp,A'',(P_{\infty},\gtb_{\infty}))}\]
is finite dimensional. It is easy to extract from Lemma \ref{cohom_for_Weyl}
that as an $\gth$-module this space is semisimple and therefore is
isomorphic to
\[\oplus_{\lambda}((M^{A}\otimes M_{\lambda,k}^{P_{0},\gtb_{0}}\otimes
M_{\lambda,k}^{P_{\infty},\gtb_{\infty}})_{\gtg(\cc_{O}^{\vee},A^{\vee})}.\]
By Lemma \ref{cohom_for_admiss}, in the last formula $\lambda$ can be chosen to
be admissible and the Verma modules can be replaced with the corresponding
admissible representations.

This proves that $\oplus_{\lambda}\Delta_{\psi^{\vee}}(M^{A}\otimes
L^{P_{0},\gtb_{0}}_{\lambda,k}\otimes
L^{P_{\infty},\gtb_{\infty}}_{\lambda,k})$ is smooth and gives a morphism

\[\Delta_{\tilde{\psi}}(M^{A})\rightarrow
\oplus_{\lambda}\Delta_{\psi^{\vee}}(M^{A}\otimes
L^{P_{0},\gtb_{0}}_{\lambda,k}\otimes
L^{P_{\infty},\gtb_{\infty}}_{\lambda,k}).\]

That this is an isomorphism can be shown in the standard way constructing the
inverse map using the formal character of $L_{\lambda,k}$, see
\cite{beil_feig_maz}.

(ii) The higher genus case is not much different. For example, pinching makes a
torus into  a sphere. Therefore in this case proof is literally the same. It
also proves  an analogue of Lemma \ref{cohom_for_admiss} for a torus. This
provides a basis for induction.

In genus $\geq 2$
at an appropriate place instead
of Proposition \ref{finofcoinv_general} one has to make reference to
Proposition
\ref{highgen_finit_afterpionching} and then use induction. $\Box$

\subsubsection{Remarks}
\label{remonquadrdeg}

(i) Meaning of Proposition \ref{statemonquadrgen} is that the dimension of the
generic
fiber of the $D-$module $\Delta_{\psi}(M^{A})$ can be calculated by the usual
combinatorial algorithm: by pinching the surface and further inserting
inserting all possible representations the problem is reduced to the case of a
sphere with three punctures and in the latter case the complete results are
available.

(ii) In the genus 0 case the  analogue of Proposition \ref{statemonquadrgen}
for generalized Weyl modules
is valid. To see this it is enough to examine  part (i) of the proof and
convince oneself that the only requirement on $M^{A}$ used there was that
$M^{A}$ be generalized Weyl module; in fact at an appropriate place instead of
Lemma \ref{cohom_for_admiss} one has to use Lemma \ref{Berma-Weyl-gener}.

(iii) Quadratic degeneration for generalized Weyl modules on the sphere allows
to write horizontal sections of the corresponding bundle as a product of vertex
operators. This will be explained in  sect.\ref{screenopcorrfuncttt}.

\section{{\bf Screening operators and correlation functions}}
\label{screenopcorrfuncttt}

In this section we will study in detail the situation described in
\ref{holondmodoncc}: we have the trivial rank 2 bundle $\ce\rightarrow\cp$, a
generalized Weyl module $M^{A}$, and
a holonomic $D-$module $\Delta(M^{A})$ on the space $\nc^{m}\times\nc^{m}$ with
fiber
$(M^{A})_{\gtg(\cp,A)}$. For the reasons which will become clear later
we replace this bundle with the dual one, its fiber being
$((M^{A})^{\ast})^{\gtg(\cp,A)}$. Denote the corresponding $D-$module by
$\Delta(M^{A})^{\ast}$. Using our results on quadratic degeneration we rewrite
horizontal sections of the corresponding bundle with flat connection as  matrix
elements of vertex operators, which serves the two-fold purpose: we find that
the differential equations satisfied by horizontal sections are provided by the
singular vectors of the corresponding Verma module and write down integral
representations for solutions to these differential equations.

\subsection{\bf Vertex operators and corelation functions}
\label{Vertex_operators_and_corelation_functions}
An alternative to the language of coinvariants in the genus zero case is the
language of {\em vertex operators}.

{\bf Definition.} A vertex operator is a $\hgtg-$morphism
\begin{equation}
\label{defvertoper}
Y:\;\cf_{\alpha \beta}^{\nc^{\ast}}\otimes V_{1}\rightarrow V_{2},
\end{equation}
where $\cf_{\alpha \beta}^{\nc^{\ast}}$ is a loop module (see
\ref{Loop_modules}) and
 $V_{1},V_{2}\in\co_{k}$ are highest weight modules. $\Box$

\bigskip

In other words, a vertex operator is an embedding
 $\cf_{\alpha \beta}^{\nc^{\ast}}\hookrightarrow
\mbox{Hom}_{\nc}V_{1}\rightarrow V_{2}$. The space
$\cf_{\alpha \beta}^{\nc^{\ast}}$ has the basis $\{F_{ij}=F_{i}\otimes
z^{j},\;i,j\in\nz\}$,
 where $\{F_{i},\;i\in \nz\}$ is a basis in $\cf_{\alpha \beta}$, see
\ref{Loop_modules}.
 Given a vertex operator $Y$, consider the generating function
\[Y(x,z)=x^{\Delta_{1}}z^{\Delta_{2}}
\sum_{i,j=-\infty}^{\infty}F_{ij}x^{-i}z^{-j},\]
the ``monodromy
coefficients'' $\Delta_{1},\Delta_{2}$ are defined by:
\[\Delta_{1}=\frac{-\lambda_{2}+\lambda_{1}+\beta}{2},\;
\Delta_{2}=\frac{-C(\lambda_{2})+C(\lambda_{1})+C(\beta)}{2},\]
where $\lambda_{i}$ is the highest weight of $V_{i}$ and
$C(\lambda)=\lambda(\lambda+2)/2$.
($\Delta_{1},\Delta_{2}$ will later appear as genuine monodromy coefficients of
a certain flat
connection.)

 The formal series
$Y(x,z)$ is, of course, an element of $Hom_{\nc}(V_{1},V_{2}\otimes
 x^{\Delta_{1}}z^{\Delta_{2}}\nc[[x^{\pm 1},z^{\pm 1}]])$. Further, for any
$g\in\gtg$ the
commutator $[g\otimes z^{n},Y(x,z)]$ is also a well-defined element of
$Hom_{\nc}(V_{1},V_{2}\otimes
 x^{\Delta_{1}}z^{\Delta_{2}}\nc[[x^{\pm 1},z^{\pm 1}]])$. For the standard
basis of $\gtg$, see

                                           \ref{defofalgebrandmodusu}, one
derives from the definition of a vertex operator that

\begin{equation}
[e\otimes z^{n},Y(x,z)]=z^{n}(-x^{2}\frac{\partial}{\partial x}+\beta x)Y(x,z),
\label{comm_vert_op_1}
\end{equation}
\begin{equation}
[f\otimes z^{n},Y(x,z)]=z^{n}\frac{\partial}{\partial x}Y(x,z),
\label{comm_vert_op_2}
\end{equation}
\begin{equation}
[h\otimes z^{n},Y(x,z)]=z^{n}(2x^{2}\frac{\partial}{\partial x}-\beta x)Y(x,z).
\label{comm_vert_op_3}
\end{equation}
We conclude that for any $g\in\gtg$ there is a differential operator $D_{g}(x)$
in $x$ such that
\begin{equation}
[g\otimes z^{n},Y(x,z)]=z^{n}D_{g}(x)Y(x,z),
\label{comm_vert_op_4}
\end{equation}

Suppose now we are given a collection of vertex operators
\[Y_{i}:\;\cf_{\lambda_{i} \mu_{i}}^{\nc^{\ast}}\otimes V_{i-1/2}\rightarrow
V_{i+1/2},\;
1\leq i\leq m.\]

The product of the corresponding generating functions
$Y_{m}(x_{m},z_{m})\cdots Y_{2}(x_{2},z_{2})Y_{1}(x_{1},z_{1})$ is a
well-defined element of
$Hom_{\nc}(V_{1/2},V_{m+1/2}\otimes
\prod_{i}x_{i}^{\Delta_{i,1}}z_{i}^{\Delta_{i,2}}
\nc[[x_{1}^{\pm 1},\ldots x_{m}^{\pm 1},z_{1}^{\pm 1},
\ldots z_{m}^{\pm 1}]])$. The matrix element
\[ <v^{\ast},Y_{m}(x_{m},z_{m}\cdots
Y_{2}(x_{2},z_{2})Y_{1}(x_{1},z_{1})v>,\;v\in V_{1/2},
v^{\ast}\in V_{m+1/2}^{\ast}\]
is, therefore, a formal Laurent series in $x_{i},z_{i},\;1\leq i\leq m$.

{\bf Definition} Suppose $Y_{i}(x_{i},z_{i}),\;1\leq i\leq m$  are as above.
Then the matrix
element
\begin{equation}
\label{defcorrfun}
\Psi(x_{1},\ldots ,x_{m},z_{1},\ldots z_{m})=
<v^{\ast},Y_{m}(x_{m},z_{m})\cdots Y_{2}(x_{2},z_{2})Y_{1}(x_{1},z_{1})v>
\end{equation}
is called {\em correlation function} if $V_{1/2},\ldots,V_{m+1/2}$ are
irreducible generalized
Weyl modules, $V_{1/2}$ is the vacuum module,
$v$ is the highest weight vector of $V_{1/2}$ and $v^{\ast}$ is the dual to the
highest weight vector of $V_{m+1/2}$. (The latter condition is meaningful in
view of the weight
space decomposition of a highest weight module.) $\Box$

\bigskip

A correlation function has been understood as a formal power series. We will
show that, in fact,
it is a holomorphic function satisfying a certain holonomic system of partial
differential
equations. In order to do that we  will interpret  vertex opeartors as
horizontal sections
of a line bundle with a flat connection provided by three modules on $\cp\times
\cp$.

\subsection {{\bf From coinvariants to vertex operator algebra}}

\subsubsection{ }
\label{coinv_corrfunsn}

We return to the setup of \ref{Thegenericlevelcase}. In the cartesian product
$\cp\times\cp$
fix coordinate system $(x,z)$. Attach to the point $x$ in the first factor the
Borel subalgebra
$\gtb_{x}$ spanned by the vectors $e_{x}=e-xh-x^{2}f,\; h_{x}=h+2xf$. This
means, in particular,
that $\gtb_{0}$ is the standard Borel subalgebra $\nc e\oplus\nc h$
(see  \ref{defofalgebrandmodusu} ) and $\gtb_{\infty}$ is the opposite one. Set
$A=
\{(0,0),(x,z),(\infty,\infty)\}$. Let $V^{A}=V^{\gtb_{0},0}_{0}\otimes
V^{\gtb_{x},z}_{1}
\otimes V^{\gtb_{\infty},\infty}_{\infty}$ be a generalized Weyl module over
$\hgtg^{A}$.
Consider the space of invariants $((V^{A})^{\ast})^{\gtg(\cp,A)}$. By Theorem
\ref{fusalggencase}
 this
space is either 0- or 1-dimensional. Suppose the latter possibility is the
case. Then by
Theorem \ref{smoothinadm} we get a line bundle with flat connection over
$\nc^{\ast}\times \nc^{\ast}$ whose fiber over the point $(x,z)\in
\nc^{\ast}\times \nc^{\ast}$ is
$((V^{A})^{\ast})^{\gtg(\cp,A)}$. There arises an embedding
\[V^{\gtb_{x},z}_{1}\hookrightarrow \mbox{Hom}_{\nc}(V^{\gtb_{0},0}_{0}\otimes
 V^{\gtb_{\infty},\infty}_{\infty},\nc).\]

The dual space $(V^{\gtb_{\infty},\infty}_{\infty})^{\ast}$ as a $\hgtg-$module
is isomorphic
to the contragredient module $(V^{\gtb_{\infty},\infty}_{\infty})^{c}$,
see \ref{defofalgebrandmodusu}. As the level is generic,
the latter module is irreducible and is, therefore, isomorphic to a certain
generalized Weyl
module
$V^{\gtb_{0},0}_{\infty}$. Hence we get an embedding
\[V^{\gtb_{x},z}_{1}\hookrightarrow \mbox{Hom}_{\nc}(V^{\gtb_{0},0}_{0},
V^{\gtb_{0},0}_{\infty}\otimes x^{\Delta_{1}}z^{\Delta_{2}}\nc[[x^{\pm 1},
z^{\pm 1}]]),\]
where $\Delta_{1},\Delta_{2}$ are monodromy coefficients of the flat
connection. We conclude that
any $w\in V^{\gtb_{x},z}_{1}$ can be looked upon as a certain generating
function
$w(x,z)=x^{\Delta_{1}-n}z^{\Delta_{2}-l}\sum_{i,j\in\nz}w_{ij}x^{-i}z^{-j}$ of
a family of
operators $\{w_{ij}\subset \mbox{Hom}_{\nc}(V^{\gtb_{0},0}_{0},
V^{\gtb_{0},0}_{\infty})$, where $(n,l)$ is a bidegree of $w$ as an element of
$V^{\gtb_{x},z}_{1}$.

\begin{lemma}
\label{vacuum_vertoper}
Suppose $v_{1}\in V^{\gtb_{x},z}_{1}$ is the highest weight vector. Then

(i) $v_{1}(x,z)$ is a generating function of a certain vertex operator as in
\ref{Vertex_operators_and_corelation_functions};

(ii) any vertex operator is obtained in this way.
\end{lemma}

{\bf Proof} is a direct and simple  calculation using definitions, see also
 \ref{prooflemma49} formula  (\ref{strofloopmodact}). $\Box$

\bigskip

Let us now relate correlation functions to horizontal sections of the bundle
built on the generalized Weyl module $M^{A}$, $\Delta(M^{A})^{\ast}$, see
beginning of sect.\ref{screenopcorrfuncttt} for notations.
Suppose that $M^{A}$ is the tensor product of ``individual'' generalized Weyl
modules \[\otimes_{i=1}^{m} V_{i}^{z_{i},\gtb_{x_{i}}}.\] Consider all possible
correlation functions
\[<v^{\ast}, v_{m}(x_{m},z_{m})\cdots v_{1}(x_{1},z_{1})v>,\]
where $v_{i}(x_{i},z_{i})$ is a generating function of a vertex operator
related
to the highest weight vector $v_{i}\in V_{i}^{z_{i},\gtb_{x_{i}}}$.
\begin{corollary}
\label{betwcorrfunandhorizsecttt}
Let $M^{A}$ be as above.
Over a suitable  open contractible subset $U$ of $\nc^{m}\times\nc^{m}$, there
is
an isomorphism between the space of horizontal sections of the  bundle
$\Delta(M^{A})^{\ast}$ and the space of correlation functions
\[<v^{\ast}, v_{m}(x_{m},z_{m})\cdots v_{1}(x_{1},z_{1})v>.\]
\end{corollary}

{\bf Proof.} Intertwining properties of vertex operators imply a correlation
function is a horizontal section of $\Delta(M^{A})^{\ast}$ in a formal sense.
This give a map in one direction. A map in the opposite direction in provided
by
quadratic degeneration, see Proposition \ref{statemonquadrgen}. $\Box$

\subsubsection{ }
\label{vertop_opervertalg}
By Lemma \ref{vacuum_vertoper} coinvariants recover vertex operators. In fact
they give us much more: the
collection of generating functions $w(x,z),\;w\in V^{\gtb_{x},z}_{1}$ affords a
kind of
{\em vertex operator algebra} structure. We will not discuss the latter in
detail
(see \cite{fren_lep_meur}) and only explain how one can get exlicit formulas
for
$w(x,z),\;w\in V^{\gtb_{x},z}_{1}$ in terms of the vertex operator $v_{1}(x,z)$
related to the
highest weight vector $v_{1}$.

For any $g\in
\gtg$ set $g(i)=g\otimes z^{i}\in\hgtg$. Define the {\em current} $g(z)$ to be
$g(z)=\sum_{i\in\nz}g(i)z^{-1-i}\in\hgtg\otimes\nc[[z^{\pm 1}]]$. Define
$g(z)^{(l)}$ to be
the $l-$th (formal) derivative of $g(z)$ with respect to $z$.
 For any  $g(z)^{(l)}$ set
\[(g(z)^{(l)})_{+}=(\frac{d}{dz})^{l}\sum_{i=0}^{\infty}g_{-i-1}z^{i},\;
(g(z)^{(l)}_{-}=g(z)^{(l)}-(g(z)^{(l)})_{+}.\]
Observe that for any $w(x,z)\in \mbox{Hom}_{\nc}(V^{\gtb_{0},0}_{0},
V^{\gtb_{0},0}_{\infty}\otimes x^{\Delta_{1}}z^{\Delta_{2}}\nc[[x^{\pm
1},z^{\pm 1}]])$
and any $g\in\gtg$, the
products $(g(z)^{(l)})_{-}w(x,z)$, $w(x,z)(g(z)^{(l)})_{+}$ are also
well-defined elements of
$\mbox{Hom}_{\nc}(V^{\gtb_{0},0}_{0},
V^{\gtb_{0},0}_{\infty}\otimes x^{\Delta_{1}}z^{\Delta_{2}}\nc[[x^{\pm
1},z^{\pm 1}]])$.

Define for any $g\in\gtg,\;w(x,z)\in \mbox{Hom}_{\nc}(V^{\gtb_{0},0}_{0},
V^{\gtb_{0},0}_{\infty}\otimes x^{\Delta_{1}}z^{\Delta_{2}}\nc[[x^{\pm
1},z^{\pm 1}]])$
\begin{equation}
\label{normalordering}
:g(z)^{(k)}w(x,z):=(g(z)^{(k)})_{-}w(x,z)+w(x,z)(g(z)^{(k)})_{+}.
\end{equation}

\begin{lemma}
\label{descendants}
 Let $g\in\gtg$, $w\in V^{\gtb_{x},z}_{1}$, $w(x,z)$ the corresponding element
of $\mbox{Hom}_{\nc}(V^{\gtb_{0},0}_{0},
V^{\gtb_{0},0}_{\infty}\otimes x^{\Delta_{1}}z^{\Delta_{2}}\nc[[x^{\pm
1},z^{\pm 1}]])$. Then

(i) $(g\cdot w)(x,z)=[g,w(x,z)]$;

(ii) $(g(-l)\cdot w)(x,z)=(1/(l-1)!):g(z)^{(l-1)}w(x,z):,\;l>0$;
\end{lemma}

{\bf Proof} is a direct calculation of matrix elements of the operator
$(g(-l)\cdot w)(x,z)$
based on the definition of the space of coinvariants. $\Box$

\subsection{ Differential equations satisfied by correlation functions}
\label{diffeqsatbycorrfuntttt}
We return to the setup of \ref{Vertex_operators_and_corelation_functions} and
consider
a correlation function
\[\Psi(x_{1},\ldots ,x_{m},z_{1},\ldots z_{m})=
<v^{\ast},Y_{m}(x_{m},z_{m})\cdots Y_{2}(x_{2},z_{2})Y_{1}(x_{1},z_{1})v>,\]
coming from the  product of vertex operators

\[Y_{i}:\;\cf_{\lambda_{i} \mu_{i}}^{\nc^{\ast}}\otimes V_{i-1/2}\rightarrow
V_{i+1/2},\;
1\leq i\leq m.\]

Using Lemma \ref{vacuum_vertoper} we assume that there are generalized Weyl
modules
$V_{i},\;1\leq i\leq m$ with highest weight vectors $v_{i},\;1\leq i\leq m$
such that
$Y_{i}(x,z)=v_{i}(x,z)$. An advantage of this point of view is that for any
collection
of elements $w_{i}\in V_{i},\;1\leq i\leq m$ we can consider the matrix element
\[<v^{\ast},w_{m}(x_{m},z_{m})\cdots w_{2}(x_{2},z_{2})w_{1}(x_{1},z_{1})v>.\]

\begin{lemma}
\label{fromprimtodesc}
For any
 $w_{i}\in V_{i},\;1\leq i\leq m$

\[<v^{\ast},w_{m}(x_{m},z_{m})\cdots w_{2}(x_{2},z_{2})w_{1}(x_{1},z_{1})v>=
D\cdot \Psi(x_{1},\ldots ,x_{m},z_{1},\ldots z_{m}),\]
where $D$ is a differential operator in $x'$s with coefficients in rational
functions in $z'$s.
\end{lemma}

{\bf Proof.} Start  with the function
\[<v^{\ast},v_{m}(x_{m},z_{m})\cdots
v_{i+1}(x_{i+1},z_{i+1})(g(-l)v_{i})(x_{i},z_{i})
v_{i-1}(x_{i-1},
z_{i-1})\cdots v_{1}(x_{1},z_{1})v>,\;l>0.\]
By Lemma \ref{descendants} (ii) it rewrites as
\[
<v^{\ast},v_{m}(x_{m},z_{m})\cdots
v_{i+1}(x_{i+1},z_{i+1})(g(z)^{(l-1)}_{-}v_{i}(x_{i},z_{i})-
v_{i}(x_{i},z_{i})g(z)^{(l-1)}_{+})
v_{i-1}(x_{i-1},
z_{i-1})\cdots v_{1}(x_{1},z_{1})v>,\;l>0.\]
Then commute all $g_{i},\; i<0$ through to the right and all $g_{i},\; i\geq 0$
to the left in
a standard way, c.f.\cite{fren_resh}  and use commutation relations
(\ref{comm_vert_op_1},\ref{comm_vert_op_2},\ref{comm_vert_op_3}). The case
$l=0$ is treated in
 a similar
and
simpler way using Lemma \ref{descendants} (i). Further argue by induction using
again
Lemma \ref{descendants}. $\Box$

\bigskip

By definition each $V_{i}$ is a quotient of a Verma module and therefore there
are elements, singular vectors in the corresponding Verma module (see
\ref{defofalgebrandmodusu})
$S_{i}\in U(\hgtg)$ such that $S_{i}v_{i}=0,\;1\leq i\leq m$. On the other
hand,
by Lemma \ref{fromprimtodesc} there are differential operators $D_{i},\;1\leq
i\leq m$ such that
\begin{eqnarray}
& &D_{i}<v^{\ast},v_{m}(x_{m},z_{m})\cdots
v_{2}(x_{2},z_{2})v_{1}(x_{1},z_{1})v>=\nonumber\\
& &<v^{\ast},v_{m}(x_{m},z_{m})\cdots
v_{i+1}(x_{i+1},z_{i+1})(S_{i}v_{i})(x_{i},z_{i})
v_{i-1}(x_{i-1},
z_{i-1})\cdots v_{1}(x_{1},z_{1})v>,\;1\leq i\leq m.\nonumber
\end{eqnarray}
We arrive to the following result.

\begin{lemma}
\label{firsthalfequat}
The correlation function
\[\Psi(x_{1},\ldots ,x_{m},z_{1},\ldots z_{m})=
<v^{\ast},v_{m}(x_{m},z_{m})\cdots v_{2}(x_{2},z_{2})v_{1}(x_{1},z_{1})v>\]
satisfies the system of equations
\begin{equation}
\label{firsthalfequat_eq}
D_{i}\Psi(x_{1},\ldots ,x_{m},z_{1},\ldots z_{m})=0,\;1\leq i\leq m.
\end{equation}
\end{lemma}

\bigskip

Observe that, although there are in general no explicit formulas for $D_{i}$,
the fact that
$[D_{i},D_{j}]=0$ is an obvious consequence of the definition.

We have obtained $m$ equations our function of $2m$ variables satisfies. The
rest is, of course,
 the Knizhnik-Zamolodchikov equations. Let us write them down explicitly.
 Recall that we can look upon
$\Psi(x_{1},\ldots ,x_{m},z_{1},\ldots z_{m})$ as a function of $z_{1},\ldots
z_{m}$ with
coefficients in a completed tensor product of $m$ $\gtg-$modules.
(The variables $x_{1}\ldots x_{m}$
are responsible for that, see
(\ref{comm_vert_op_1},\ref{comm_vert_op_2},\ref{comm_vert_op_3}).)
 For any $A=\sum_{s}a_{i}\otimes b_{s}
\in\gtg\otimes\gtg$ denote by $A_{ij},1\leq i,j\leq m$
an operator  acting on the $m-$fold tensor product of $\gtg$-modules by the
formula
\[A_{ij}\cdot w_{1}\otimes\cdots w_{m}=\sum_{s}w_{1}\otimes
a_{s}w_{i}\otimes\cdots b_{s}w_{j}
\otimes\cdots w_{m}.\]
The formula (\ref{comm_vert_op_4}) implies that $A_{ij}$ is a differential
operator
 in $x_{i},x_{j}$.
Set $\Omega=ef+fe+h^{2}/2$.

\begin{lemma} (\cite{knizh_zam})
\label{knzameq}

The correlation function
$\Psi=\Psi(x_{1},\ldots ,x_{m},z_{1},\ldots z_{m})$ satisfies the  system of
Knizhnik-Zamolodchikov
 equations
\begin{equation}
\label{knzameq_eq}
(k+2)\frac{\partial}{\partial z_{i}}\Psi=\sum_{j\neq
i}\frac{\Omega_{ij}}{z_{i}-z_{j}}\Psi,\;
1\leq i\leq m.
\end{equation}
\end{lemma}

There is no need to prove this lemma here as one can repeat
word for word the known proofs.
 However we point out that if one considers a highest weight module $V$ as a
module
over the semi-direct product of $\hgtg$ and the Virasoro algebra $Vir$ then $V$
is annihilated
by the element $d/dz-L_{-1}$, where $L_{-1}$ is one of the Sugawara elements.
One then shows that
  the singular vectors
$(d/dz-L_{-1})v_{i}$, where $v$ is a highest weight vector of $V_{i}$, give
rise to
the equations (\ref{knzameq_eq}) in exactly the same way the singular vectors
$S_{i}$ gave
rise to the equations (\ref{firsthalfequat_eq}). An immediate consequence of
this proof is
that the system of equations (\ref{knzameq_eq},\ref{firsthalfequat_eq}) is
consistent.

\subsection{Screening operators and integral representations of correlation
functions}
\label{scr_and_int_repr}
\subsubsection{ }
\label{motivation}
Suppose a  function $\Psi_{old}=\Psi(x_{1},\ldots ,x_{m},z_{1},\ldots z_{m})$
 is the matrix element of the product
of vertex operators
\[\Psi_{old}=<v^{\ast},\pi\circ Y(x_{m},z_{m})\cdots
Y_{1}(x_{1},z_{1})v_{o}>,\]

\[Y_{i}:\;\cf_{\lambda_{i} \mu_{i}}^{\nc^{\ast}}\otimes V_{i-1/2}\rightarrow
V_{i+1/2},\;
1\leq i\leq m,\]
satisfying the same conditions as the expression  in (\ref{defcorrfun}),
see\ref{Vertex_operators_and_corelation_functions}, except
that instead of assuming that $V_{m+1/2}$ is a generalized Weyl module we
assume that
$V_{m+1/2}$ is a contragredient Verma module, see \ref{defofalgebrandmodusu}.
(Why ``old'' will become clear in a moment.) It is easy to see that
 $\Psi_{old}=\Psi(x_{1},\ldots ,x_{m},z_{1},\ldots z_{m})$ satisfies the same
system of equations
(\ref{firsthalfequat_eq},\ref{knzameq_eq}). Suppose in addition that there is a
projection
 $\pi:\; V_{m+1/2}\rightarrow W$ onto another contragredient Verma module $W$.
Denoting
by $w^{\ast}$ an element dual to the highest weight vector $w\in W$ one can
consider the
matrix element
\[\Psi_{new}=<w^{\ast},\pi\circ Y(x_{m},z_{m})\cdots
Y_{1}(x_{1},z_{1})v_{0}>.\]
We again observe that $\Psi_{old}$ is a solution to the same system
(\ref{firsthalfequat_eq},\ref{knzameq_eq}). This new solution can be calculated
 as follows.

There arises the dual map $\pi^{\ast}:\; W^{\ast}\rightarrow V_{m+1/2}^{\ast}$
and by definition
   there is an element $S$ of $U(\hgtg_{>})$
such that $\pi^{\ast}(w^{\ast})= Sv^{\ast}$. We now take the definition of
$\Psi_{new}$, replace
in it $\pi^{\ast}(w^{\ast})$ with $Sv^{\ast}$ and get

\begin{equation}
\label{replachwvect}
\Psi_{new}=<S\cdot v^{\ast},\circ Y(x_{m},z_{m})\cdots
Y_{1}(x_{1},z_{1})v_{0}>.
\end{equation}

 Then we  commute $S$ through to the right.
 The intertwining properties of vertex operators tell us
that
\begin{equation}
\label{fromoldtonew}
\Psi_{new}=S^{t}\cdot\Psi_{old},
\end{equation}
where $^{t}$ signifies the canonical
antiinvolution an a Lie algebra
($g_{1}g_{2}\cdots g_{n}\rightarrow g_{n}g_{n-1}\cdots g_{1}$)
and the action  is determined by the following condition: if $g\in\gtg$ then
\[(g\otimes
z^{n})\cdot\Psi_{old}=\sum_{i=1}^{m}D_{g}(x_{i})z_{i}^{n}\Psi_{old},\]
see (\ref{comm_vert_op_4}).

We intend to use (\ref{fromoldtonew}) in the case when $\pi$ and therefore $S$
do not exist!

\subsubsection{Screening operators }
Let $V_{\lambda_{\infty},k}$ be a highest weight module and
$v\in V_{\lambda_{\infty},k}$ a highest weight vector.
 If the obvious integrality conditions are satisfied
then the vectors $f^{\lambda_{\infty}+1}v,\;( e\otimes
z^{-1})^{k-\lambda_{\infty}+1}v$ are
singular and give rise to embeddings of the type $W\hookrightarrow
V_{\lambda_{\infty},k}$.
Now take 3 highest weight modules $V_{\lambda_{i},k},\; i=0,1,\infty$ attach
them
to 3 point in $\cp$ and consider the space of coinvariants
\[(\otimes_{i=0,1,\infty}V_{\lambda_{i},k}^{\gtb_{i},i})_{\gtg(\cp,\{0,1,\infty\})}.\]
Of course an embedding $W\hookrightarrow V_{\lambda_{\infty},k}$ gives rise to
a map
\[(W^{\gtb_{\infty},\infty}\otimes_{i=0,1}V_{\lambda_{i},k}^{\gtb_{i},i})_{\gtg(\cp,\{0,1,\infty\})}
\hookrightarrow
(\otimes_{i=0,1,\infty}V_{\lambda_{i},k}^{\gtb_{i},i})_{\gtg(\cp,\{0,1,\infty\})}.\]
It is remarkable that even if the embedding $W\hookrightarrow
V_{\lambda_{\infty},k}$ does not
exist the last map still does. In the language of vertex operators this
phenomenon was explained
in great detail in \cite{feig_mal}.

Therefore with each of the formal singular vectors --
$f^{\lambda_{\infty}+1}v\mbox{ or }( e\otimes z^{-1})^{k-\lambda_{\infty}+1}v$
-- we have
associated an operator acting on coinvariants. Call these operators {\em
screenings} and denote
them $R_{1}$ and $R_{0}$ respectively.

Let us calculate the action of the  screenings explicitly. By definition
\newline
$R_{j}(\otimes_{i=0,1,\infty}V_{\lambda_{i},k}^{\gtb_{i},i})_{\gtg(\cp,\{0,1,\infty\})}$ only
depends on $V_{\lambda_{\infty},k}$ so we will be simply writing
$R_{j}(V_{\lambda_{\infty},k})$.
Now formulas for the related singular vectors
 ($f^{\lambda_{\infty}+1}v,\;( e\otimes z^{-1})^{k-\lambda_{\infty}+1}v$) and a
very simple
calculation using the formulas (\ref{param_eq_line_1},\ref{param_eq_line_2})
give the following
result:
\begin{eqnarray}
\label{act_scr_vert_op_1}
R_{1}((V_{m}^{0},V_{n}))=(V_{m-1}^{1},V_{n})\\
\label{act_scr_vert_op_2}
R_{1}((V_{m}^{1},V_{n}))=(V_{m+1}^{0},V_{n})\\
\label{act_scr_vert_op_3}
R_{0}((V_{m}^{0},V_{n}))=(V_{m+1}^{1},V_{n})\\
\label{act_scr_vert_op_4}
R_{0}((V_{m}^{1},V_{n}))=(V_{m-1}^{0},V_{n})
\end{eqnarray}

Suppose we are given 2 generalized Weyl modules and a vertex operator acting
between them.
Suppose in addition that this vertex operator is related to a highest weight in
the third
generalized Weyl module, say $(V_{m}^{\epsilon},V_{n})$.
Theorem \ref{fusalggencase} tells us that given such a vertex operator
our screenings give us all the others of the type $(V_{i}^{\alpha},V_{n})$ --
we cannot only
change the value of $n$.
But then there is the standard screening operator -- $S$ -- which takes care of
$n$, see
e.g. \cite{feig_fr_1}. So these three -- $R_{1},R_{2},S$ -- provide us with all
vertex operators.
 This has
an important application to the calculation of correlation functions.

Start with a simple correlation function given by the product of vertex
operators, each of which
is characterized by the condition $m=0$. Then applying $S$ an appropriate
number
 of times one gets all
vertex operators and, hence,
all correlation functions in spirit of Varchenko-Schechtman, see\cite{awtsyam}.

Now take a Varchenko-Schechtman
correlation function $\Psi_{old}$. It comes from a product of vertex operators:
\[\Psi_{old}=<v^{\ast}, Y(x_{m},z_{m})\cdots Y_{1}(x_{1},z_{1})v_{o}>,\]

\[Y_{i}:\;\cf_{\lambda_{i} \mu_{i}}^{\nc^{\ast}}\otimes V_{i-1/2}\rightarrow
V_{i+1/2},\;
1\leq i\leq m.\]

Let $W_{i},\;0\leq i\leq m$, be words on 2 letters $R_{1}$ and $R_{2}$.
Replacing $V_{i+1/2}$
with $W_{i}(V_{i+1/2})$ we get a new correlation function $\Psi_{new}$. Doing
this with all
$\Psi_{old}$ and sufficiently many $W_{i},\;0\leq i\leq m$ we get all solutions
to
(\ref{firsthalfequat_eq},\ref{knzameq_eq}). In principle all these solutions
can be written
down explicitly. It is especially simple to do so in the case when we keep
$V_{i+1/2},\;0\leq i\leq m-1,$ and only change $V_{m+1/2}$.

So assume that $\Psi_{old}$ is as above and replace $V_{m+1/2}$ with
$R_{j}(V_{m+1/2}),\; j=0,1$.
Then by (\ref{replachwvect}) one is to expect that
\[\Psi_{new}=<X^{\alpha}\cdot v^{\ast},\circ Y(x_{m},z_{m})\cdots
Y_{1}(x_{1},z_{1})v_{0}>,\]
where $X$ is either $e$ or $f\otimes z$ if $j=1$ or 0 resp., and $\alpha$ is
either $\lambda+2$
or $k-\lambda+2$ resp., where $\lambda$ is
 the highest weight of $V_{m+1/2}$.

Of course if $\alpha$ is not a nonnegative integer  then the last formula does
not make much sense.
Nevertheless using it and (\ref{fromoldtonew}) as a motivation we arrive to
\[\Psi_{new}=X^{\alpha}\Psi_{old}.\]
Now the left-hand side of the last equality does make sense:
$X$ is a first order differential operator,
see \ref{motivation}, therefore  we can set in a rather straightforward manner
\[X^{\alpha}\Psi_{old}=\int t^{-\alpha-1}\{\exp(-Xt)\Psi_{old}\}dt\]
and get a nice integral operator, for details see \cite{feig_mal}.

This procedure can be easily iterated to provide the functions
\begin{equation}
\label{our_solut}
\int \prod_{i=1}^{n}t_{i}^{-\alpha_{i}-1}\{\exp(-X_{1}t_{1})\exp(-X_{2}t_{2})
\cdots \exp(-X_{n}t_{n})\Psi_{old}\}\prod_{i=1}^{n}dt{i},
\end{equation}
where $X_{1},X_{2},...$ is either $e,f\otimes z,e,...$ or $f\otimes
z,e,f\otimes z,...$.

\begin{lemma}
\label{our_sol_stae}
Functions (\ref{our_solut}) are solutions to
(\ref{firsthalfequat_eq},\ref{knzameq_eq}).
\end{lemma}

{\bf Proof} is  same as the proof of the analogous statement in
\cite{feig_mal1}. In fact it is an easy exrcise to make the heuristic arguments
which have lead us to the formula (\ref{our_solut}) into a precise proof.
$\Box$

\bigskip

Integrating functions (\ref{our_solut}) with respect to $x'$s
(or doing something similar but more esoteric) one is supposed to get the
Dotsenko-Fateev correlation functions for the Virasoro algebra. It would be
interesting to do
this explicitly and compare the result with the calculations in \cite{f_g_p_p}.

\begin{conjecture}

(i) Formulas (\ref{our_solut}) provide all solutions to the system
(\ref{firsthalfequat_eq},\ref{knzameq_eq}).

(ii) If the level $k$ is rational, then there arises a subbundle of  the bundle
in question, the one with fiber $((L^{A})^{\ast})^{\gtg(\cp,A}$, where
$L^{A}$ is the corresponding admissible representation. We conjecture that in
this case formulas (\ref{our_solut}) actually give horizontal sections of the
latter bundle.

 \end{conjecture}


\begin{thebibliography}{99}

\bibitem{awata}Awata H.,Yamada Y., {\em Fusion rules for the Fractional Level
$\widehat{\gtsl(2)}$ Algebra}, KEK-TH-316 KEK Preprint 91-209, January 1992

\bibitem{awtsyam} Awata H., Tsuchia A., Yamada Y., Nucl.Phys., {\bf
B365}(1991), 680-698

\bibitem{ba_soch} Bauer M., Sochen N. Singular vectors by fusions in
$A_{1}^{(1)}$ Saclay preprint SPht/91-117 (1991)

\bibitem{beil_feig_maz} Beilinson A., Feigin B., Mazur B. {\em Introduction to
algebraic field
theory on curves}, preprint

\bibitem{bern_beil} Beilinson A., Bernstein J. {\em Localisation de
$\gtg-$modules},
C.R.Acad.Sci.Paris, {\bf 292} (1981), 15-18

\bibitem{bpz} Belavin A.A., Polyakov A.M., Zamolodchikov A.B., Nuclear Physics
{\bf B241} (1984) 333-380


\bibitem{bryl_kash} Brylinski J-P., Kashiwara M. {\em Kazhdan-Lusztig
conjecture and holonomic systems}, Invent.Math., {\bf 64}
(1981), 387-410

\bibitem{drinf} Drinfeld V.G. Proc.Int.Congr.Math. Berkeley  {\bf 1} 1986
\bibitem{deodgabbkac} Deodhar V.V., Gabber O., Kac V.G., Adv.in Math {\bf 45}
(1982) 92-116
\bibitem{fink} Finkelberg M. Fusion categories, Ph.D. thesis, Harvard
university,1993

\bibitem{feig_fuchs} Feigin B., Fuchs D. Journal of Geom. and Phys.
{\bf 5} (1988) 209-235

\bibitem{feig_fr} Feigin B., Frenkel E. Phys.Lett.B {\bf 246} (1990) 75
\bibitem{feig_fr_1} Feigin B., Frenkel E., in Physics and Mathematics of
Strings, V.G.Knizhnik
Memorial Volume, eds. L.Brink, et.al, 271-316, World Scientific, Singapore 1990

\bibitem{fei_fuchs} Feigin B.L., Fuchs D.B.  in Representations of Lie groups
and related topics, eds. A.M.Vershik, D.P.Zhelobenko, 465-554, Gordon and
Breach, New York 1990
\bibitem{feig_fuchs_2} Feigin B.L., Fuchs D.B. J.Geom.Phys. {\bf 5}(1988) n.2

\bibitem{feld} Felder G. Nucl.Phys. {\bf 317} (1989)

\bibitem{fkw} Frenkel E., Kac V., Wakimoto M. Comm.Math.Phys. {\bf 147}(1992)
295-328


\bibitem{feig_mal1} Feigin B., Malikov F. Letters in Math.Phys. {\bf 31} 1994,
315-325
\bibitem{feig_mal} Feigin B., Malikov F., Advances in Sov.Math. {\bf 17}(1993)
15-63


\bibitem{fren_resh} Frenkel I.B., Reshetikhin N.Yu Comm.Math.Phys. {\bf
146}(1992) 1-60
\bibitem{fren_lep_meur} Frenkel I., Lepowski J., Meurman A.{\em Vertex Operator
 Algebra and the Monster}, Academic Press, Inc 1988
\bibitem{fuchs} Fuchs D.B. Funkc. Anal. i Ego Pril. {\bf 23} (1989) 2, 81-83

\bibitem{f_g_p_p} Furlan P., Ganchev A.Ch., Paunov R., Petkova V.B.,
 Phys.Letters {\bf 267} (1991)
63; {\em Solutions of
 the Knizhnik-Zamolodchikov equations with rational isospins and the reduction
to the minimal models}, preprint CERN-TH.6289/91,
 accepted for publication in {\em Nucl.Phys.B}

\bibitem{hitch} Hitchin N. Duke Math.J., {\bf 54} (1987), n.1
\bibitem{kac_kazhd} Kac V.G., Kazhdan D.A., Adv.in Math.
{\bf 34} (1979) 97-108
\bibitem {ioh_mal} Iohara K., Malikov F., Modern Phys.Lett.A {\bf 8}(1993),
No.38 3613-3624

\bibitem{jimbo} Jimbo M. Lett.Math.Phys. {\bf 10} 63-69 (1985)

\bibitem{kac_wak} Kac V.G., Wakimoto M. Proc. Nat'l Acad. Sci. USA
{\bf 1988} 4956

\bibitem{kazh_luszt} Kazhdan D., Lusztig G. Duke Math.J. {\bf 62} 21-29
\bibitem{knizh_zam} Knizhnik V.G., Zamolodchikov A.B. Nucl.Phys. {\bf B 247}
(1984) 83-103

\bibitem{kul_resh} Kulish P.P., Reshetikhin N.Yu. Lett.Math.Phys. {\bf 18}
143-149 (1989)


\bibitem{malff}Malikov F.G., Feigin B.L., Fuchs D.B.,
Funkc.Anal.i ego Pril. {\bf 20}(1988) 2, 25-37

\bibitem{mal}Malikov F.,
  Infinite
Analysis - Proceedings of the RIMS Research Project 1991 Part B,
 623 - 645, World Scientific Co.
Pte. Ltd.

\bibitem{mal_2} Malikov F. Leningrad Math.Journal {\bf 2} (1991) 269-286

\bibitem{mat_walt}  Mathieu P., Walton M. Prog. Theor. Phys. Suppl.
{\bf 102} (1990) 229

\bibitem{moorseib} Moore G., Seiberg N. Comm.Math.Phys. {\bf 123} (1989)
177-254

\bibitem{peters} Petersen J.L., Rasmussen J., Yu M. Conformal blocks for
admissible representations in $\gtsl(2)$ current algebra, NBI-HE-95-16,
hep-th/9504127

\bibitem{pr_seg} Pressley A., Segal G. Loop groups, Oxford Mathematical
Monographs, 1988

\bibitem{ramg} Ramgoolam S. New modular Hopf algebras related to rational $k$
$\widehat{\gtsl}_{2}$ YCTP-P2-93, hep-th/9301121

\bibitem {resh_tur} Reshetikhin N.Yu., Turaev V.G. Inv.Math. {\bf 103} 547-597
(1991)

\bibitem{seg} Segal G., in Proceedings of XI International Congress on
Mathematical
Physics, 17 - 27  July, 1988 Swansea, Adam Hilger Bristol and New York, 1988
\bibitem{sch_varch} Schechtman V., Varchenko A., Invent.Math. {\bf 106}(1991),
139-194

\bibitem{ts_u_ya} Tsuchiya A., Ueno K., Yamada Y.,
Adv. Studies in Pure Math. {\bf 19} (1989) 459-566
\bibitem{verl} Verlinde E. Nucl.Phys. {\bf B300} (1988) 360

 \end{thebibliography}
\end{document}